\documentclass{article}

\usepackage[UKenglish]{babel}
\usepackage{graphicx}
\usepackage{fullpage}
\begin{document}
\title{Numerical approximations to extremal metrics on toric surfaces}
\author{R. S. Bunch and S. K. Donaldson}
\date{\today}
\maketitle
%\begin{document}
   \newtheorem{thm}{Theorem}
\newtheorem{lem}{Lemma}
\newtheorem{cor}{Corollary}
\newtheorem{prop}{Proposition}
\newtheorem{conj}{Conjecture}
\newcommand{\bC}{{\bf C}}
\newcommand{\bR}{{\bf R}}
\newcommand{\cH}{{\cal H}}
\newcommand{\cP}{{\cal P}}
\newcommand{\bP}{{\bf P}}
\newcommand{\oa}{\overline{a}}
\newcommand{\ua}{\underline{a}}
\newcommand{\dbd}{\overline{\partial}}
\newcommand{\cA}{{\cal A}}
\newcommand{\SDiff}{{\rm SDiff}}
\newcommand{\cJ}{{\cal J}}
\newcommand{\cJint}{{\cal J}_{{\rm int}}}
\newcommand{\Riem}{{\rm Riem}}
\newcommand{\grad}{{\rm grad}}
\newcommand{\Ricci}{{\rm Ric}}
\newcommand{\ut}{\underline{t}}
\newcommand{\ux}{\underline{x}}
\newcommand{\utheta}{\underline{\theta}}
\newcommand{\bZ}{{\bf Z}}
\newcommand{\Xcx}{X_{{\rm cx.}}}
\newcommand{\Xsymp}{X_{{\rm symp}}}
\newcommand{\Xalg}{X_{{\rm alg}}}
\newcommand{\oP}{\overline{P}}
\newcommand{\Vol}{{\rm Vol}}
\newcommand{\Av}{{\rm Av}}
\newcommand{\cL}{{\cal L}}
\newcommand{\Lie}{{\rm Lie}}
\newcommand{\ueta}{\underline{\eta}}
\newcommand{\uepsilon}{\underline{\epsilon}}

\section{Introduction}
This is a report on numerical work on  {\it extremal Kahler metrics} on toric
surfaces, using the ideas developed in \cite{kn:D1}\footnote{Although the
simpler results described here were obtained by the second author in 2004, before those of \cite{kn:D1}.}.  These are special
Riemannian metrics on manifolds of $4$ real, or $2$ complex dimensions.  They can be regarded as  solutions
of a fourth order partial differential equation for a convex function on a polygon. We refer to the preceding article \cite{kn:D2} for further background. Motivation
for this kind of  work can be found in at least three different sources:
\begin{itemize}
\item  Great effort is expended on the proof of abstract existence theorems
for such metrics. Only very rarely does one expect to find explicit formulae
for the solutions. It is interesting to explore the geometry of the solutions
using numerical methods.
\item It is conceivable that specific numerical results of this general nature
could be relevant to physics \cite{kn:Doug}.
\item The numerical work illuminates the analytical techniques and difficulties
involved in the abstract existence theory. This is true both of the approximate
numerical solutions and of the approximation scheme itself, which is founded
on a link between algebraic geometry and asymptotic analysis.
\end{itemize}

The results described here are obtained using a computer programme which
can be applied in principle to any polygon and in practice works effectively
for a reasonable range of cases. We  focus on four such cases: particular
$p$-sided
polygons with $p=5,6,7,8$, which correspond to  complex surfaces
constructed by \lq\lq blowing up'' the projective plane  $p-2$ times.
The status of the rigorous existence theory is different in the four cases:
\begin{itemize}
\item For the hexagon, the metric we seek is a Kahler-Einstein metric whose
existence was proved by Tian, Yau  and Siu in the 1980's. Using a rather different approach,
 Doran {\it et al} \cite {kn:DHHKW} also recently found numerical approximations
 to this metric.
\item For the pentagon the metric we seek was proved to exist very recently
by Chen, LeBrun and Weber \cite{kn:CBW}
\item For the octagon, the metric we seek has constant scalar curvature and
the abstract existence follows from  recent work of the second author.
\item For the heptagon, the metric we seek is extremal but not constant scalar
curvature, and as far as we know is not at present covered by any general
existence theorem.
\end{itemize}

But we emphasise that these cases are just a small sample of the ones that
can be treated successfully, and  there are many other cases of the fourth
kind, where the numerical results outrun, at present, the abstract existence theory.

In each of these four cases we give a selection of geometrical results: again
these form a small sample of what one can do.
The pentagon is particularly interesting because, as Chen, Le Brun and Weber
showed, it yields a new (conformally-Kahler) solution to the Riemannian Einstein
equations. (More precisely, the extremal metric we approximate is very close to  another
extremal metric which is conformal to an Einstein metric.) We verify this
numerically. from which one could go on to obtain information about the Einstein metric.

In the final section, we make some comments about ways in which these techniques
might be taken further in the future.

\

We are grateful to Julien Keller for many discussions about this work, and
much crucial technical assistance. In particular, the C++ version of the
programme to find geodesics was written by Keller.

\section{The set-up}

\subsection{Algebraic metrics}
Our starting point is an integral Delzant polygon $P$, contained in the square
of side length $k$. So for example, this is a picture of our hexagon, with
$k=6$.

\

\

\

\

\begin{picture}(150,150)(-70,-20)
\put(0,0){\line(1,0){75}}
\put(0,0){\line(0,1){75}}
\put(75,0){\line(1,1){75}}
\put(0,75){\line(1,1){75}}
\put(75,150){\line(1,0){75}}
\put(150,75){\line(0,1){75}}
\put(5,30){{\bf -1}}
\put(30,5){{\bf -1}}
\put(38,95){{\bf -1}}
\put(95,38){{\bf -1}}
\put(137,105){{\bf -1}}
\put(105,137){{\bf -1}}
\put(-7,-9){(0,0)}
\put(152,152){(6,6)}
\put(75,-9){(3,0)}
\end{picture}

When we write, for example, \lq\lq the hexagon with $k=12$'', we mean the
same figure scaled up by a factor of $2$. So, for the hexagon, we can take
$k$ to be any even number (because we need to work with polygons having integral
vertices). This use of language differs slightly from that in \cite{kn:D2},
but we hope the meaning will be clear. We also indicate in the diagram the self-intersection
numbers of the $2$-spheres in the manifold corresponding to the edges; in
this case all $-1$.

\

Next we consider a collection of positive numbers $(a_{\nu})$ indexed by
the lattice points $\nu$ in the closed polygon $\oP$. Then we define a \lq\lq
Kahler potential''
 $$  \phi(\ut) = \log \left( \sum_{\nu} a_{\nu} e^{\nu.\ut}\right). $$
 This a convex function (of a variable $\ut\in \bR^{2}$) whose derivative maps
 onto the polygon $P$. The Kahler potential defines an \lq\lq algebraic metric''
 on the toric surface $X$ corresponding to $P$. The Legendre transform of $\phi$ is a convex function $u$ on
 $P$. The scalar curvature $S$ is given by the expression $S(u)=-\sum_{ij} \frac{u^{ij}}{\partial x_{i}\partial x_{j}}$. If $u$ were to satisfy the partial differential equation
 \begin{equation} S(u)=-A,
 \end{equation} for an affine linear function $A$, then our metric would
 be extremal. This occurs in two very special cases:
 \begin{itemize}
 \item If $P$ is the $m_{1}\times m_{2} $ rectangle, with origin as one vertex,
 then we can take the product of binomial coefficients
 $$  a_{\nu_{1}, \nu_{2}}= \left(\begin{array}{c}m_{1}\\ \nu_{1}\end{array}\right)
 \left(\begin{array}{c}m_{2}\\ \nu_{2}\end{array}\right), $$
 for integers $0\leq \nu_{1}\leq m_{1}, 0\leq \nu_{2}\leq m_{2}$.
 \item If $P$ is the triangle with vertices $(0,0), (k,0), (0,k)$ we can
 take  the multinomial co-efficents
 $$   a_{\nu_{1}, \nu_{2}}= \frac{k!}{\nu_{1}! \nu_{2}! (k-\nu_{1}-\nu_{2})!}.
 $$
 \end{itemize}
 The metrics we describe are, respectively, the product of round metrics on  $S^{2}\times
 S^{2}$ and the Fubini-Study metric on $\bC\bP^{2}$. These are extremal metrics
 but of course we have not really gained anything new. In general, we cannot
 hope to find an extremal metric {\it exactly}, in this way, but we can hope
 to find good {\it approximations}, and we can hope that the approximations
 improve when we increase $k$. For example, this is a collection of co-efficients
 $a_{\nu}$ which give a good approximation to the extremal metric in the
 case of the hexagon, for $k=6$:

 \begin{picture}(500,240)(-70,-20)
 \put(0,0){\line(1,0){105}}
\put(0,0){\line(0,1){100}}
\put(105,0){\line(1,1){100}}
\put(0,100){\line(1,1){100}}
\put(100,200){\line(1,0){105}}
\put(205,100){\line(0,1){100}}
 \put(5,5){2.06}
 \put(35,5){7.17}
 \put(65,5){7.17}
 \put(95,5){2.06}
 \put(5,35){7.17}
 \put(35,35){38.9}
 \put(65,35){63.5}
 \put(95,35){38.9}
 \put(125,35){7.17}
 \put(5,65){7.17}
 \put(35,65){63.5}
 \put(65,65){163.6}
 \put(95,65){163.6}
 \put(125,65){63.5}
 \put(155,65){7.17}
 \put(5,95){2.06}
 \put(35,95){38.9}
 \put(65,95){163.6}
 \put(95,95){257.2}
 \put(125,95){163.6}
 \put(155,95){38.9}
 \put(185,95){2.06}
 \put(35,125){7.17}
 \put(65,125){63.5}
 \put(95,125){163.6}
 \put(125,125){163.6}
 \put(155,125){63.5}
 \put(185,125){7.17}
 \put(65,155){7.17}
 \put(95,155){38.9}
 \put(125,155){63.5}
 \put(155,155){38.9}
 \put(185,155){7.17}
 \put(95,185){2.06}
 \put(125,185){7.17}
 \put(155,185){7.17}
 \put(185,185){2.06}
\end{picture}

\

\

In this case there is a symmetry under a dihedral group of order 12, so we
can reduce to only 6 independent coefficients. For clarity we have just
written these to a few significant figures. When $k$
is larger the data $(a_{\nu})$ is best kept confined to a computer file,
and we will just record the maximum and minimum of the $a_{\nu}$. It may
be helpful to recall here from \cite{kn:D3} that the numbers $(-\log a_{\nu})$ can
be thought of, roughly, as a \lq\lq discrete approximation'' to the symplectic
potential $u$. Saying much the same thing in a different way, if we view
 $u$ as a function on $P$ depending on parameters $a_{\nu}$, then at a
given point $\ux \in P$ the value of $u$ is strongly influenced by the $a_{\nu}$
for $\nu$ close to $\ux$, but only little by those for far away $\nu$. These somewhat
vague statements become more precisely true when the parameter $k$ is large.

\

If $\alpha$ is any affine-linear function then changing the coefficients $a_{\nu}$ to
$ e^{\alpha(\nu)} a_{\nu}$ does not change the metric described by the data.
We can normalise the $a_{\nu}$ to take account of this, as follows. Let $(p_{1},
p_{2})$ be the centre of mass of the polygon. We arrange
that \begin{equation}\sum_{\nu} \log a_{\nu}=0\ , \ \sum_{\nu} (\nu_{1}-p_{1}) \log a_{\nu}=0\ , \ \sum_{\nu} (\nu_{2}-p_{2})
\log a_{\nu}=0.\end{equation}

\subsection{Decomposition of the curvature tensor}
 
Here we want to extend slightly the discussion of the curvature tensor of
a toric manifold in \cite{kn:D2}, for the two-dimensional case.

 We begin in a more general setting and consider a complex Kahler surface $Z$ with an anti-holomorphic
isometric involution $\sigma$ fixing a submanifold $\Sigma$ of two real dimensions.
Thus $\Sigma$ is a \lq\lq real form'' of $Z$. We want to understand the curvature
tensor of $Z$ at a point $p$ of $\Sigma$. We have $TZ_{p}=T\Sigma_{p}\oplus N_{p}$ where $N$ is the
normal bundle. Fix an orientation of $\Sigma$, at least locally. This determines
an orientation of the normal bundle so, using the metric, we can define complex structures on
each of $TS$ and $N$, making them complex line bundles. With suitable choice
of signs the resulting complex structure $J_{-}$ on $TX_{p}$ induces the
{\it opposite} orientation on $Z$ from that of the Kahler structure. Thus
we have a unit anti-self-dual $2$-form $\theta_{-}$ at $p$ and we can write
$$   \Lambda^{-}_{p}= \bR \theta_{-} \oplus L, $$
say, where $L$ is a two-dimensional oriented real vector space. Thinking
of
$TZ_{p}$ as the complexification of $TS_{p}$ we see that there is a natural
map from the trace-free symmetric $2$-tensors $s^{2}_{0}=s^{2}_{0}T^{*}_{p}$
to the trace-free Hermitian forms $\Lambda^{1,1}_{0}$ on $Z$ at $p$. In a
standard  basis this is just the obvious map which takes a symmetric real matrix
to a Hermitian matrix. It is easy to check that the image of this is $L$,
so we can identify $L$ with $s^{2}_{0}$. 

 Now we use
the well-known decomposition of the curvature tensors of Riemannian $4$-manifolds,
and in particular Kahler surfaces. In the general Riemannian case, the curvature tensor, viewed as a symmetric element of
$\Lambda^{2}\otimes \Lambda^{2}$, decomposes into four pieces according to
the $\pm$-self-dual splitting of the $2$-forms, and two of these are equivalent
by the symmetry. The $\Lambda^{+}\otimes \Lambda^{-}$
block is equivalent to the trace-free Ricci tensor. When the $4$-manifold
is Kahler the $\Lambda^{+}\otimes \Lambda^{+}$ piece (the self-dual Weyl
tensor $W^{+}$) is entirely determined
by the scalar curvature. The scalar curvature also determines the trace of
the piece $W^{-}$ in $\Lambda^{-}\otimes \Lambda^{-}$. Further, since the Ricci tensor has type $(1,1)$, the
$\Lambda^{+}\otimes\Lambda^{-}$ piece lies in $ \langle \omega\rangle\otimes \Lambda^{-}$, where
$\omega$ is the Kahler form (which is self-dual).

Putting these observations together, we see that the curvature tensor of
$Z$ at $p$ has at most the following  six irreducible components
\begin{enumerate}
\item The scalar curvature $S$.
\item A component of the trace-free Ricci tensor in $\langle\omega\rangle\otimes L\subset
\Lambda^{+}\otimes \Lambda^{-}$. This can be regarded as an element $\rho$ of $s^{2}_{0}$.
\item A component of the trace-free Ricci tensor in the $1$-dimensional space
$\langle \omega \rangle\otimes
\langle \theta_{-}\rangle \subset \Lambda^{+}\otimes \Lambda^{-}$.
\item A component of the anti-self-dual Weyl tensor $W_{-}$ in $\langle \theta_{-}\rangle
\otimes \langle \theta_{-}\rangle\subset \Lambda^{-}\otimes \Lambda^{-}$.
\item A component of $W^{-}$ in $\langle \theta_{-}\rangle \otimes L\subset
\Lambda^{+}\otimes \Lambda^{+}$.
\item A component $w$ of $W^{-}$ in the 2-dimensional space $s^{2}_{0}(s^{2}_{0})$.
\end{enumerate}

It is an exercise to check that, since the curvature tensor is preserved
by the action of $\sigma$ on $TZ_{p}$, which reverse the Kahler form, the components (3) and (5) vanish.
We can identify the component $\rho$ in (2), which determines the trace-free Ricci
tensor, with a quadratic differential on $\Sigma$--- an element of $T^{*}S\otimes_{\bC}T^{*}S$.
Likewise, we can identify the component $w$ in (6) with a quartic differential on
$\Sigma$. The scalar component (4) can be written as
$  S- 6 K$ where $K$ is the Gauss curvature of the induced metric on $\Sigma$.
The conclusion is that the curvature tensor is determined by $S,\rho, K,
w$.

\

We now apply this to our toric surface $X$, taking $\Sigma$ to be the real form
$X_{\bR}$. We work in the local co-ordinates $t_{1}, t_{2}$,  but the point of the
above discussion is to see that the tensors we obtain extend smoothly over
this compact surface. (If we want to use the complex description, with quadratic
and quartic differentials, we should go to the oriented cover.)
The four-dimensional Riemannian metric associated to the convex function
$\phi(t_{1}, t_{2})$ is given
by
$$   \sum \phi_{ab} dt_{a} dt_{b} + \phi_{ab} d\theta_{a} d\theta_{b}. $$
Here   $\theta_{1}, \theta_{2}$ are angular
coordinates and the $\phi_{ab}$ are the second derivatives
of $\phi$ with respect to the $t$-variables. The $\theta_{a}$
coordinates never appear in the calculations and, because of the torus action,
the curvature tensor over the whole of $X$ is obviously determined by its
restriction to $X_{\bR}$.

We define a tensor
\begin{equation}  F_{ab cd}=\phi_{a b c d}- \sum_{i j} \phi_{ab i} \phi_{cd j} \phi^{ij},\end{equation}
where $(\phi^{ij})$ is the inverse of the Hessian $\nabla^{2} \phi=(\phi_{ab})$. 
The tensor $F$ is symmetric in the pairs of indices $(a b)$ and $(cd)$ and also under
interchange of the two pairs. The scalar curvature is
$$  S= \sum_{a b c d} F_{ab cd} \phi^{ab} \phi^{cd}. $$ 
The quadratic
differential $\rho$ corresponds to the trace-free Ricci tensor and
is given by  
$$  \rho_{ab}=\sum_{c d} F_{ab cd} \phi^{cd}-  \frac{S}{2} \phi_{ab}. $$
The other scalar component is determined by the Gauss curvature $K$ of the
metric restricted to the $(t_{1}, t_{2})$ variables, which is
$$ K= (F_{11 22}- F_{1212}) \det(\phi_{ab})  $$
The final component $w$, in the space of quartic differentials, is obtained by
projecting $F$ to the  the kernel  of these contractions.  We will 
call this the \lq\lq Weyl component'', although strictly the Weyl curvature
of the $4$-manifold also contains contributions from $S$ and $K$. With our
conventions we have
\begin{equation}  \vert {\rm Riem} \vert^{2} = \frac{1}{3} S^{2} + \frac{1}{24} (S-6K)^{2}
+ \vert w\vert^{2} + \frac{1}{2} \vert \rho \vert^{2}. \end{equation}
{\it Examples}
\begin{itemize}
\item If $X$ is $\bC\bP^{2}$, with the standard Fubini-Study metric scaled to $S=1$, then \begin{equation} 
S=1 \ \ \rho=0 \ \ \ \vert w\vert^{2} =0 \ \ \ \ K=1/6\end{equation}
and $\vert {\rm Riem}\vert^{2}= 1/3$.
\item If $X$ is $S^{2}\times S^{2}$, with each factor having the round metric
of the same area then, scaled to $S=1$,
\begin{equation} S=1\ \ \ \rho= 0 \ \ \ \vert w\vert^{2}= 1/8 \
\ \ \ K=0 \end{equation}
and $\vert {\rm Riem}\vert^{2}= 1/3+ 1/24+1/8=1/2$.
\end{itemize}

We note also that Chern-Weil theory gives a global relation
\begin{equation}
\int_{X} \vert {\rm Riem} \vert^{2} - S^{2}     = (2\pi)^{2} (p-6),\end{equation}
where $p$ is the number of vertices of the polygon.

\

We will also consider the \lq\lq Bach tensor''. According to Chen, LeBrun
and Weber,  \cite{kn:CBW}, this
can be written, on any Kahler surface, as the trace-free part of
$$ S \ \rm{Ricci} + i\partial \overline{\partial} S, $$
where $S$ is the scalar curvature and ${\rm Ricci}$ is the Ricci form. If
the metric is extremal this is a harmonic $(1,1)$ form. In our setting we
can write the Bach tensor as a quadratic differential given
by 
$$  B_{ab}= S \rho_{ab} + S_{ab} - (S_{ij} \phi^{ij}) \phi_{ab}, $$
where $S_{ab}$ denotes the partial derivative $\frac{\partial^{2} S}{\partial t_{a}\partial t_{b}}$. If the scalar curvature $S$ is positive and we consider
the conformally related metric $\hat{g}= S^{-1} g$ then the trace-free Ricci
tensor of $\hat{g}$ is  $S B_{ab}$. In particular if the Bach tensor vanishes
then this conformal metric is Einstein. All this is explained in \cite{kn:CBW}.

Suppose we have any trace-free symmetric $2$-tensor $\theta_{ab}$ on $X_{\bR}$. This defines a $(1,1)$
form $\theta_{ab} d\tau^{a} d\overline{\tau}^{b}$ on $X$. The condition that
this form is harmonic is simply that 
$$  \frac{\partial \theta_{ab}}{\partial t_{c}}= \frac{\partial \theta_{ac}}{\partial
t_{b}}. $$
 This equation is similar to that
defining the {\it holomorphic} quadratic differentials on $X_{\bR}$, in terms of
the conformal structure on $X_{\bR}$.  More precisely,
the equation can be written as $\overline{\partial} \theta + L(\theta)=0$, where
$L$ is a certain algebraic operator, determined by $\phi$. Let us call the
solutions of this equation \lq\lq $\phi$-holomorphic quadratic differentials''. Then
we see that for an extremal metric the Bach tensor can be viewed as a $\phi$-holomorphic
quadratic differential on $X_{\bR}$. One can show further that, if the Bach
tensor vanishes, then the product $S w$ is a \lq\lq $\phi$-holomorphic quartic differential'' on $X_{\bR}$, in a similar sense. Just as for genuine holomorphic
differentials, the zeros of $\phi$-holomorphic differentials occur with positive
signs, so the number of zeros, counted with multiplicity, is fixed by topology.

       \subsection{Integration}
       
       The foundation of our approximation scheme is the numerical integration
       of functions on the manifold $X$. Since the functions will all be
       invariant under the torus action the integration is really in $2$ dimensions
       rather than $4$, but the presence of the fixed points of the action, corresponding
       to the boundary of the polygon, means that this reduction is not
       completely straightforward. Suppose given a coefficent set $(a_{\nu})$.
       The functions we  want to integrate
       are given explicitly, by standard formulae, in the $t$ coordinates,
       so one approach would be to fix a lattice covering a large square
       in the $(t_{1}, t_{2})$-plane and approximate the integrals by the
       corresponding sums. There are two drawbacks to this. First, we will
       get errors from the contributions to the integrals outside the square.
       Second, in these $\ut$ co-ordinates the calculations will
       involve ratios with very small numerator and denominator when $\vert
       \ut\vert$ is large, which leads to numerical errors. Of course this
       is a reflection of the fact that these co-ordinates are not valid
       over the whole of $X$. 
       
       To get around these difficulties we proceed as follows. Suppose we
       have a fixed convex function $U_{0}$ on the polygon $P$, computable in
       the $\ux$-coordinates. The derivative of $U_{0}$ gives a map $DU_{0}:P\rightarrow
       \bR^{2}$. On the other hand the derivative of the Kahler potential
       $\phi$, defined by the $(a_{\nu})$ gives a map $D\phi: \bR^{2} \rightarrow
       P$. So we have a composite map $D\phi\circ DU_{0}: P\rightarrow P$. Let $J(\ux)$ be the  determinant of the derivative of this map. Fix a small lattice spacing $h=k/N$ and let $\ux^{(\alpha)}$
       run over the points in the intersection $P\cap N^{-1}\bZ^{2}$. Let
       $f$ be a $T^{2}$-invariant function on $X$, regarded as a function
       $f(\ut)$ on $\bR^{2}$.
       The integral of $f$ over $X$ is approximated by the sum
       \begin{equation} h^{2}  \sum_{\alpha} J(\ux^{(\alpha)}) \ f(DU_{0}(\ux^{(\alpha)})). \end{equation}
       If we think of $f$ as a function on $P$, then we are approximating
       the integral by a weighted average of the values at the array of points
       $\tilde{\ux}^{(\alpha)}=D\phi\circ DU_{0} (\ux^{(\alpha)})$. Recall
       that we have a system of charts covering $X$, each chart associated
       to a vertex of $P$.
       For each point $x^{(\alpha)}$ we choose a suitable  chart
       containing the point $DU_{0}(\ux^{(\alpha)})$. Thus,  after a linear transformation, we can suppose       that this vertex is the origin and that $P$ coincides, near the origin,
      with the standard model $\{x_{i}>0\}$. The corresponding chart is formed
      by taking complex co-ordinates 
      $$  z_{a} = \exp(\frac{1}{2}( t_{a}+ i\theta_{a})), $$
      and the chart is chosen so that we work at a point with $\vert z_{a}\vert
      \leq 1$.
      In these coordinates
      
      $$ \phi= \sum_{\nu} a_{\nu_{1}, \nu_{2}}  \vert z_{1} \vert^{2\nu_{1}}
      \vert z_{2} \vert^{2\nu_{2}}. $$
      We perform all our calculations in these coordinates to evaluate the
      contribution to the sum (8) from the point $\ux^{(\alpha)}$. 
       
We choose the function $U_{0}$ to be $b u_{0}$ where $u_{0}$ is Guillemin's \lq\lq admissible
symplectic potential'' 
$$  u_{0}= \sum_{r} ( \lambda_{r}- c_{r}) \log (\lambda_{r}- c_{r}), $$
in the notation of \cite{kn:D2}, and $b$ is a positive constant. If we took
$b=1$ the pattern of points $\tilde{\ux}^{(\alpha)}$ would be similar to
the original pattern $\ux^{(\alpha)}$. Indeed, if it happened that $\phi$
were the Legendre transform of $u_{0}$ then these sets would be identical. With this choice,  the integration scheme would have errors of order $N^{-1}$ arising
from the boundary of $P$.  So we take $b>1$, which  has the effect
of increasing the density of the set $\tilde{\ux}^{(\alpha)}$ near the boundary
and reducing
the error term. In fact we have taken $b=2$ which seems to give sufficient
accuracy. (We  test the scheme by calculating
integrals whose value we know exactly, such as the total volume of the manifold
or Chern-Weil integrals: for a reasonable number of integration points we
typically get errors   in the range $10^{-6}- 10^{-4}$.) The whole
procedure is illustrated in Figure 1. The hexagon on the left shows the
lattice of points $\ux^{(\alpha)}$ in $P$. Those that map into a particular chart, corresponding to the top left vertex, appear in the shaded region in
the left hand diagram. This chart is represented by the small diagram in
the middle, with co-ordinates in  the square corresponding to $\vert z_{1}\vert,
\vert z_{2}\vert$. The points where we calculate, in this chart, appear as
dots in the square. The hexagon on the right displays the set of points
$\tilde{\ux}^{(\alpha)}$, with those corresponding to the given chart shaded.

\begin{figure}[!h]
  \caption{The integration scheme}
  \centering
    \includegraphics[scale=0.4]{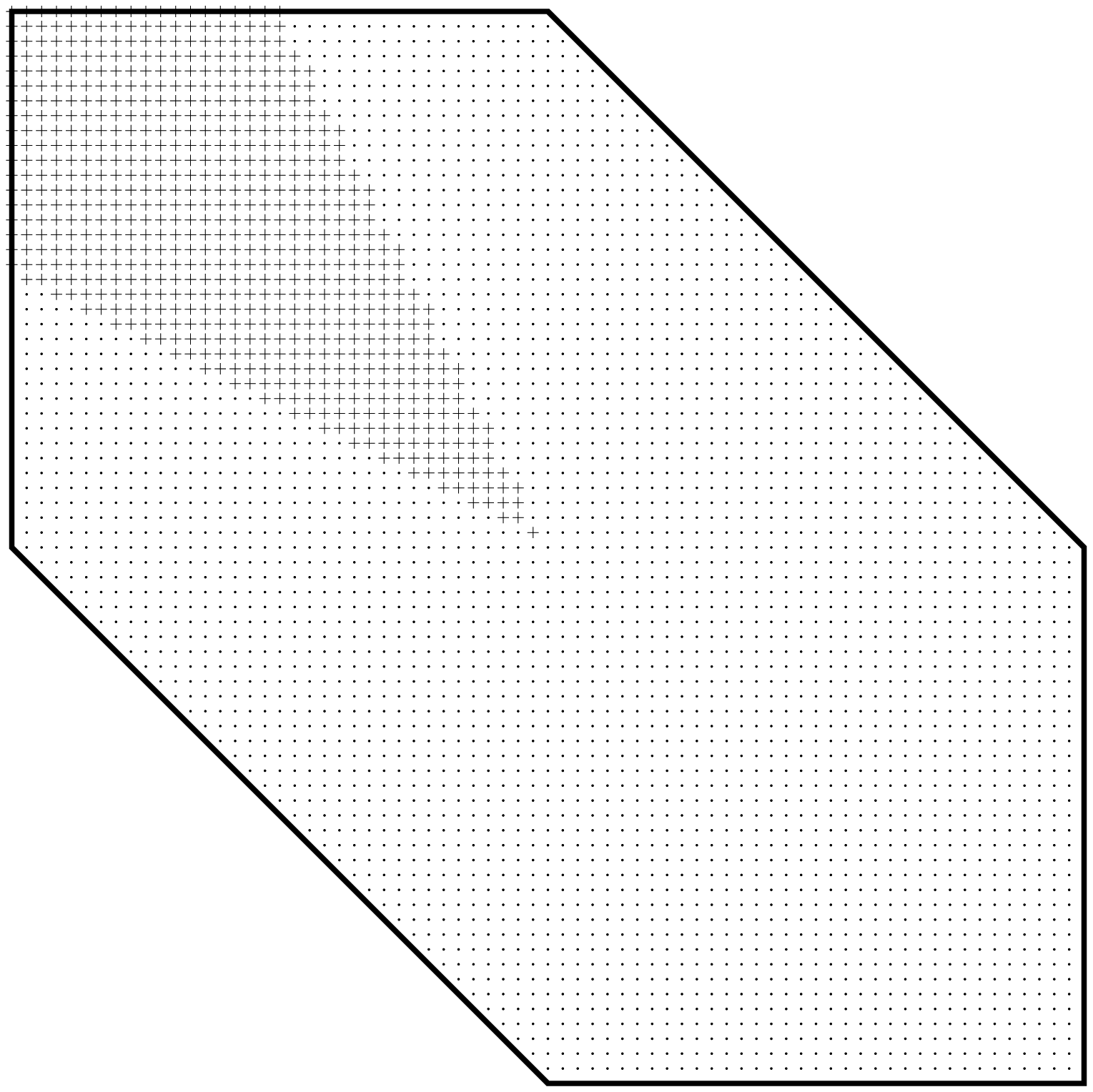}
    \includegraphics[scale=0.15]{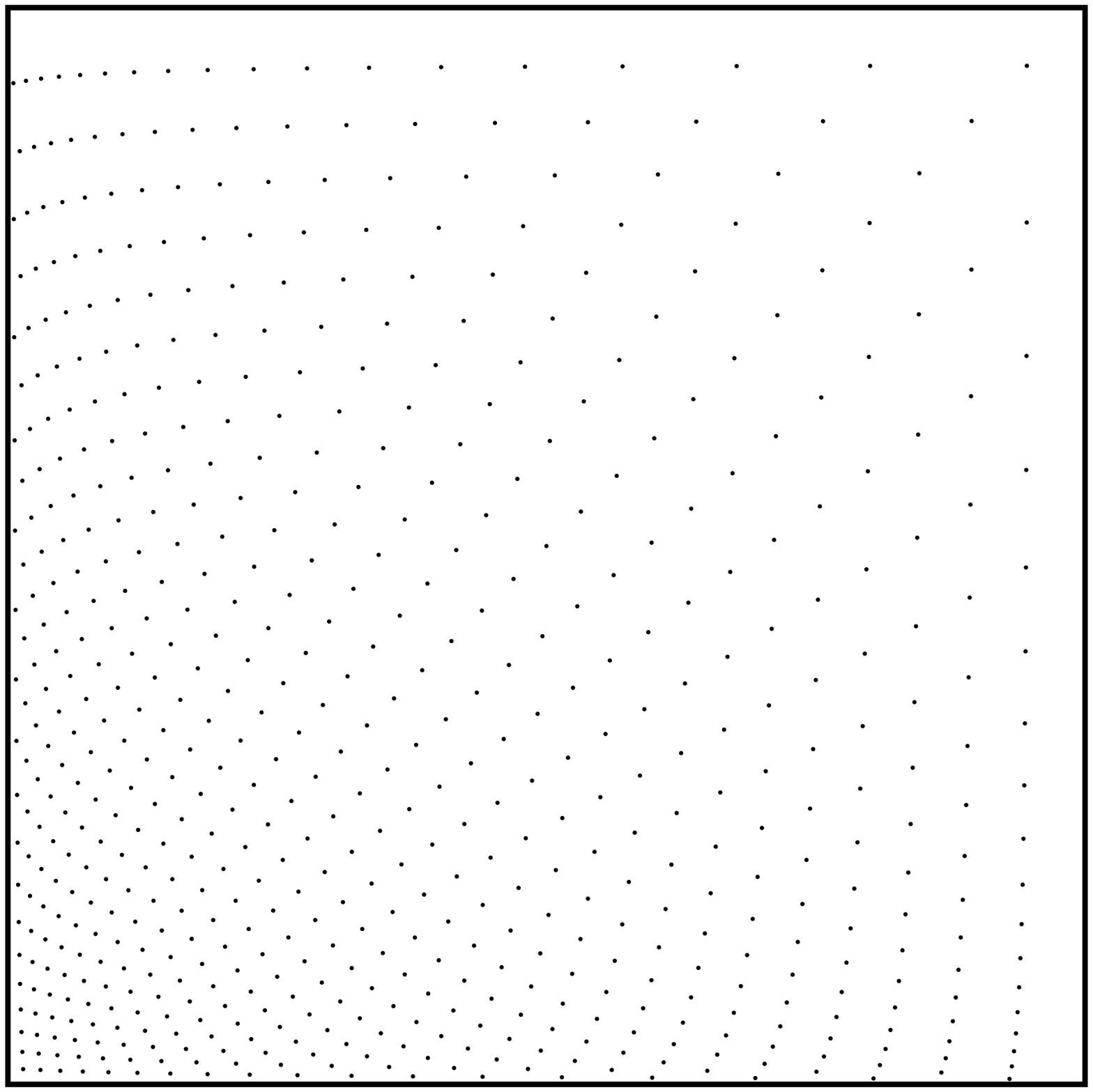}
    \includegraphics[scale=0.4]{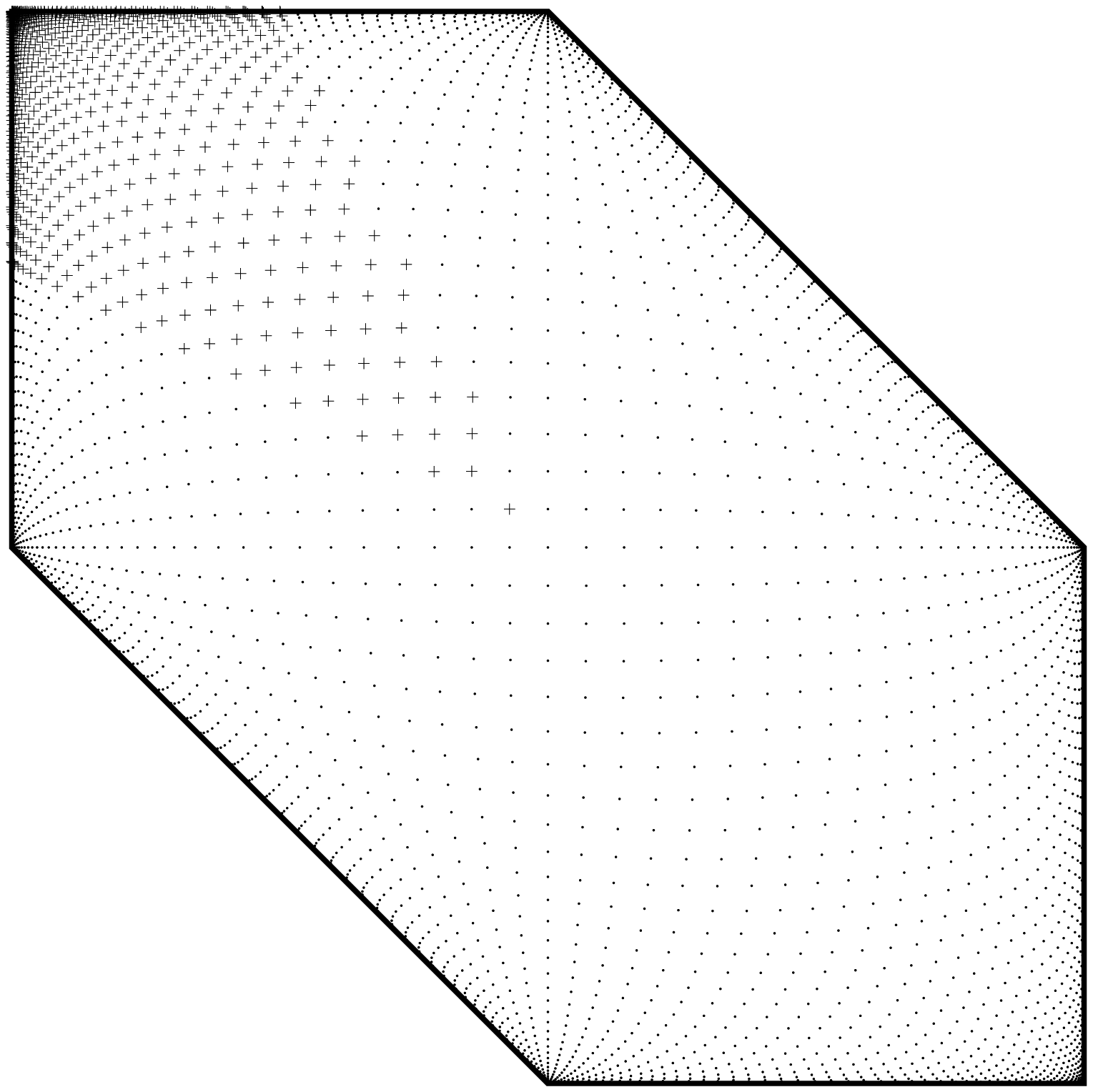}
\end{figure}

\section{Numerical algorithms: balanced metrics and refined approximations}

Our problem is to find coefficient sets $(a_{\nu})$ such that the corresponding
metric is close to an extremal metric. We consider two different conditions
on   the
$(a_{\nu})$ which we can hope have {\it exact solutions}, close to the extremal
metric. The first is the   \lq\lq balanced'' condition. The general
context for this involves the $L^{2}$ norm on the spaces of holomorphic sections
of line bundles, but in this toric situation we can write everything out explicitly, ignoring the algebro-geometric
interpretation if we please. Given $(a_{\nu})$ we write
$$ D(\ut) = \sum_{\nu} a_{\nu} e^{\nu. \ut}, $$
do $\phi(\ut)=\log D(\ut)$. For each lattice point $\mu\in P$ we have a function
$$    f_{\mu}(\ut)= \frac{ e^{\mu.\ut}}{D(t)}. $$
Then we perform the integrals
$$  I_{\mu}= \int_{\bR^{2}} f_{\mu}(\ut) \ \det \nabla^{2} \phi \ d\ut . $$
We
have written these as integrals over $\bR^{2}$ but of course we evaluate
them numerically using the scheme described above. Now
set \begin{equation}a^{*}_{\nu}= I_{\nu}^{-1}.\end{equation} We
say that the coefficient set $(a_{\nu})$, or the metric it defines, is {\it balanced} if $(a^{*}_{\nu})$ is equivalent to $(a_{\nu})$ in the sense that, for some affine linear function
$\alpha$,
$$   a^{*}_{\nu}= e^{\alpha(\nu)} a_{\nu}. $$
General theory\footnote{Strictly, some of this theory has not been written
down for the case at hand, but to simplify exposition we will ignore this
here.} tells us that the solutions of this balanced condition are essentially unique,
and that if an extremal metric exists then balanced metrics exist for large
enough $k$ and converge to the extremal metric as $k\rightarrow \infty$.
Moreover  there is a simple algorithm for finding these balanced metrics.
Recall  that we fixed a way (2) to normalise the coefficients to take account
of the action of the affine linear functions. Let $\cH$ denote the set of
normalised coefficients.
We take any set of normalised co-efficients $(a_{\nu})$ and compute $a^{*}_{\nu}$. Then
we normalise the $a^{*}_{\nu}$ to get a new set $(a^{**}_{\nu})$. This defines
a map $T:\cH\rightarrow \cH$ and the balanced metrics are exactly the fixed
points of this map. These can be found by starting with any point $\ua= (a_{\nu})$ in $\cH$
and taking the limit of iterates $T^{n}(\ua)$ as $n\rightarrow \infty$.

This procedure provides a robust and straightforward way to find balanced metrics numerically,
and these give reasonable approximations to the extremal metrics.
The drawback is that they do not give very accurate approximations, for practical
values of $k$. One can only
expect the balanced metric to differ by $O(k^{-1})$ from the extremal metric.
Thus, as in \cite{kn:D1} Sec. 2.2.1,  we seek to extend the ideas to \lq\lq refined approximations''. Recall that, on the polygon $P$,  the extremal metric is a solution of the equation $S=A$ where
$A$ is an affine-linear function of the symplectic coordinates $x_{i}$, determined by the elementary
geometry of the polygon. The extremal metric is a minimum of the modified Mabuchi functional ${\cal F}$ which
is a functional on the space of all Kahler metrics. By restricting to the
algebraic metrics we get a function $F$ on $\cH$. We say that a coefficient
set $(a_{\nu})$ is a {\it refined approximation} if it is an absolute minimiser
of $F$. 

The condition that $(a_{\nu})$ is a critical point of $F$ can be expressed
directly, without explicit reference to the Mabuchi functional. We regard
the scalar curvature $S$ as a function on $\bR^{2}$, in the obvious way,
and write $\tilde{A}$ for the composite $A\circ D\phi$ so, expressed on $\bR^{2}$,
the extremal equation is $S=\tilde{A}$. For each $\mu$ we form the integrals
\begin{equation}   \eta_{\mu}= \int_{\bR^{2}} (S-\tilde{A}) f_{\mu}(\ut) \ \det \nabla^{2} \phi
\ d\ut. \end{equation}
Then $(a_{\nu})$ is a critical point for $F$ if and only if all the $\eta_{\mu}$
vanish. 

The abstract theory of these refined approximations is less clear-cut than
for the balanced metrics. But it is reasonable to expect that, assuming an
extremal metric exists,
\begin{itemize}\item When $k$ is sufficiently large there is a unique refined
approximation and it is close to the balanced metric.
\item As $k\rightarrow \infty$ the refined approximation converges to the
extremal metric faster than any power $k^{-m}$.
\end{itemize}
We refer to the discussion, for a related problem, in \cite{kn:D1}, Sections
4,5 and move on to the
main issue,  for our present
purposes, which is finding these
refined approximations in practice. (Notice that the \lq\lq refined approximation''
is an exact concept, and of course we can only hope to find approximations
to it.)

We will give a very sketchy discussion to motivate the procedure for finding
refined approximations which we
use. Imagine we start with the balanced metric and write $\hat{S}_{0}$ for
the difference $S_{0}-A$, where $S_{0}$ is the scalar curvature of the balanced metric. We compute the \lq\lq error coefficients'' $\ueta= (\eta_{\nu})$ as above.
Suppose we have a collection of small numbers $(\epsilon_{\nu})$, which we
will refer to as \lq\lq correction coefficients'',  labelled
by the lattice points in $\oP$. Then we can modify the map $T$ to a new map
$T_{\uepsilon}$ by changing the definition (9) to
$$   a^{*}_{\nu} = (1+\epsilon_{\nu}) I_{\nu}^{-1}. $$
If $\uepsilon$ is small we can hope that iteration of $T_{\uepsilon}$ will
converge to a fixed point, close to the balanced metric. This fixed point
has some scalar curvature $S_{1}$ and we set $\hat{S}_{1}= S_{1}-A$. One
can argue that $\hat{S}_{1}$ is approximately $\hat{S}_{0} - \sum_{\nu} \epsilon_{\nu}
f_{\nu}$. Now let the $L^{2}$-projection of $\hat{S}_{0}$ onto the span of the
functions $f_{\nu}$ be $\sum \sigma_{\nu} f_{\nu}$. If we choose $\epsilon_{\nu}=\sigma_{\nu}$
then we can hope that $\hat{S}_{1}$ is small compared with $\hat{S}_{0}$,
{\it i.e} that the fixed point of $T_{\uepsilon}$ is a better approximation
to the extremal metric than the balanced metric. But by definition,
$$  \sum_{\nu} M_{\mu \nu} \sigma_{\nu} = \eta_{\mu}, $$
where
$$  M_{\mu \nu}= \int_{\bR^{2}} f_{\mu}(\ut) f_{\nu}(\ut) \det(\nabla^{2} \phi) \ d\ut.
$$
So, if we can invert the matrix $M$ we should define the correction coefficients
$\uepsilon$ in terms of the error coefficients $\ueta$ by
$$   \uepsilon=M^{-1} \ueta, $$
and iterate $T_{\uepsilon}$ to find a fixed point which is a better approximation
than the balanced metric. 

Now the matrix $M$ will be very large and, even when practical, inverting it exactly will probably
not be very wise.  As explained in \cite{kn:D1}, Sect.4.2, 
      the matrix $M$ can be expected to be an approximation to the heat operator
      operator
      $e^{-\Delta/k}$ on the manifold. So we expect that there are many eigenvalues
      of $M$ which are slightly less than $1$, but also some very small eigenvalues.
      (This  is confirmed by numerical tests.) So the exact inverse
      of $M$ will be very large. It is better to invert only on the span
      of the eigenvalues which are close to $1$.   But on this span the inverse is close to $M$ itself.
     The upshot is that, rather than inverting the matrix, we could just set $\ueta=
     c \uepsilon $ for some appropriate fixed constant $c$. (In fact, after
     experimentation,  we set $c=0.75$.) Then we can still hope that the
     fixed point of $T_{\uepsilon}$ is a better approximation than the balanced
     metric. 
     
     With this motivation in place we describe the actual procedure we use. At each stage we have a set of \lq\lq correction
     co-efficients'' $\epsilon_{\nu}$. At the outset these are set to $0$.
     We iterate the map $T_{\uepsilon}$ until we are close to a fixed
     point. ( Here \lq\lq closeness'' is judged by the criterion that the ratio
     $a^{**}_{\nu}/ a_{\nu}$ differs from $1$ by at most some small number
     $c'$ which, after experimentation, we set to $.0008$.) Then we calculate the \lq\lq error co-efficients'' $\eta_{\nu}$
     by performing the integrals (10). We modify the correction coefficients
     to a new set
     $$  \tilde{\epsilon}_{\nu}= \epsilon_{\nu} + .75 \eta_{\nu}, $$
     and proceed to iterate $T_{\tilde{\uepsilon}}$, and so on. A fixed point
     of the whole procedure should correspond to the refined approximation,
     with $\eta_{\nu}=0$.

     In practice, leaving aside the very sketchy rationale for the procedure,         this 
     process does  converge, for practical purposes, in all the
      examples we have tested. 
      After a very long while, we find experimentally that the procedure    typically         oscillates
      between two values, but these differ only in the seventh decimal digit
      so effectively we have reached a fixed point.
      
      On the negative side,  this procedure is painfully
      slow when $k$ is large. Just as in \cite{kn:D1},% (where a similar
    %   procedure was used to approximate Calabi-Yau metrics),
     we can hope that the technique
  could be  improved  with a more detailed and sophisticated
      analysis.% (perhaps
     % including a less crude approximate inversion of the matrix $M$).

      We graph the $L^{2}$ norm of $S-A$, plotted
      on a logarithmic scale, 
       under
      this procedure; for two of the  cases discussed in the next Section. %The graphs for
 %     these two cases (and all others  we have tried) are similar in
 %     nature. 
 %The $L^{2}$ norm always decreases. 
 There is a short initial
      period in which the decrease is very rapid. Then there is long intermediate
      period where  the decrease is approximately
      exponential,  until the norm is quite close  to its limiting value. After this the decrease is much slower. (We expect that the first period
is that in which nonlinear effects are important and that in the second period
the iteration is close to the iteration of a linear map, which in turn should
be close to the linearisation of the \lq\lq Calabi flow''.)
\begin{figure}[!h]
 % \caption{The $L^{2}$ norm of the error in the iteration: the pentagon %with $k=20$.}
  \centering
    \includegraphics[scale=0.4]{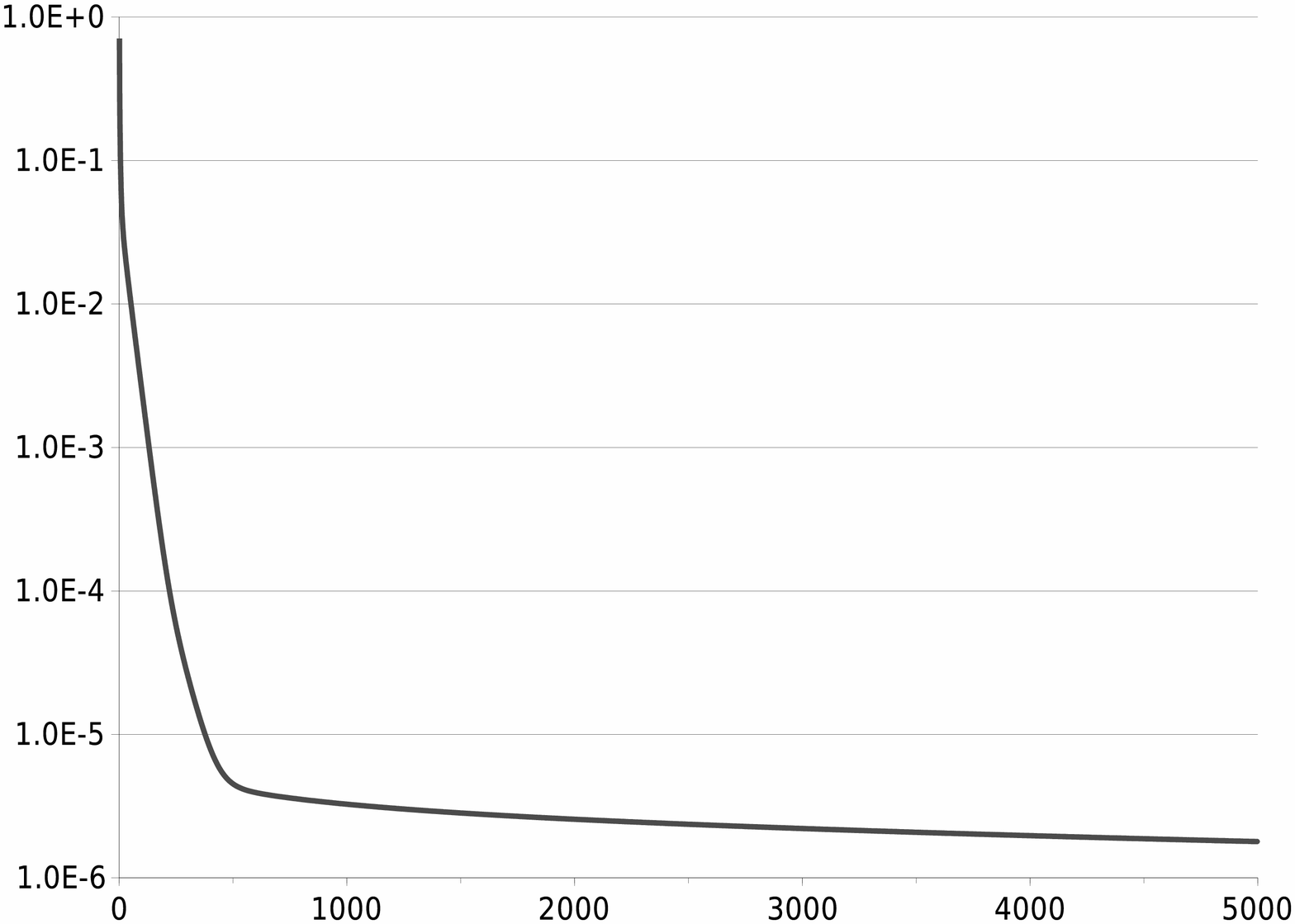}
\end{figure}

\begin{figure}[!h]
  \caption{Evolution of $L^{2}$ error: the pentagon with $k=20$ (above) and
  the heptagon with $k=45$ (below)}
  \centering
    \includegraphics[scale=0.4]{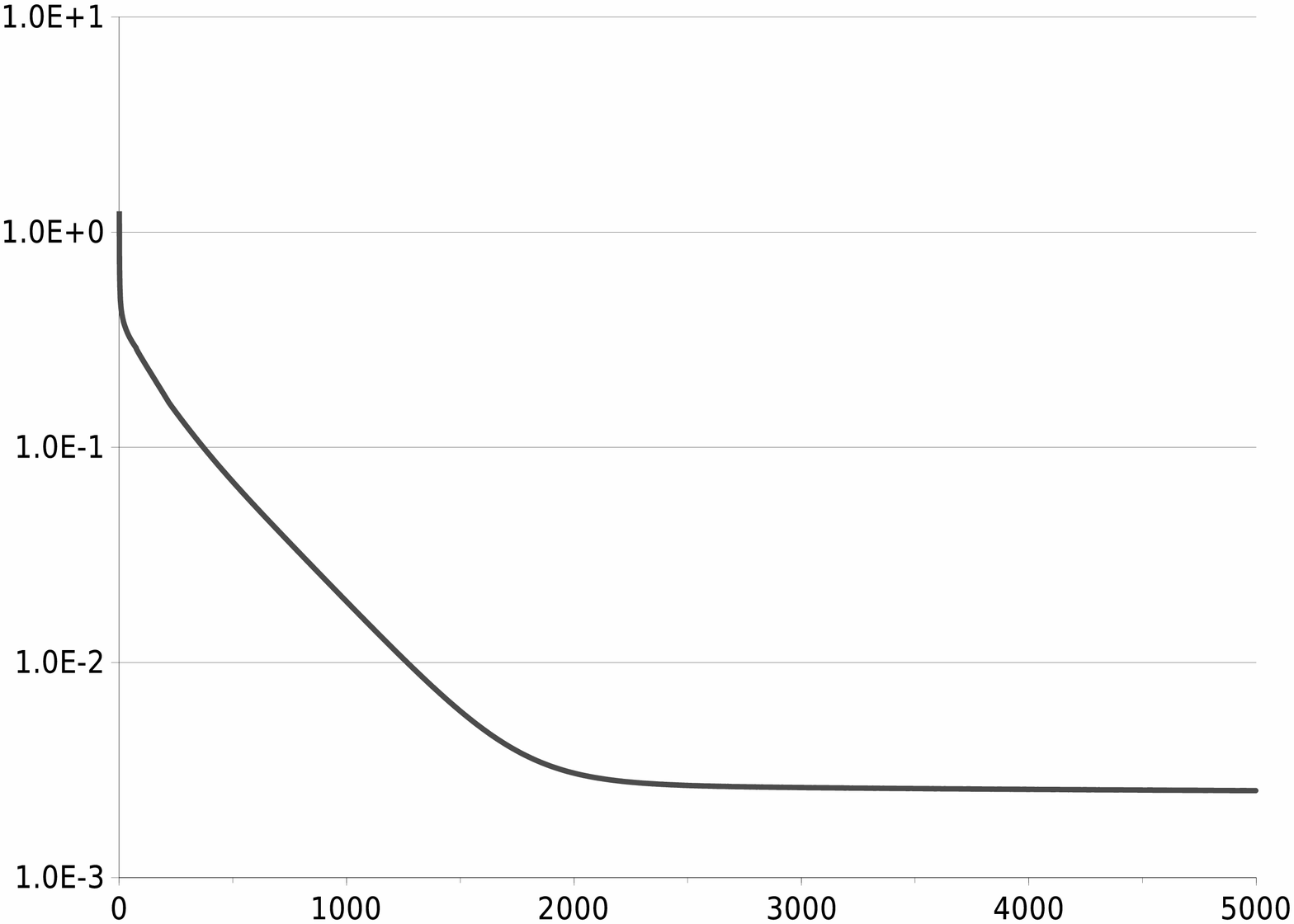}
\end{figure}
       
     \clearpage
       
   \section{Numerical results}

   Suppose we have  an algebraic metric which we hope to be a
   good approximation to an extremal metric. To compare across different
   values of $k$ we rescale so that the average value of $S$ is $1$. Setting
   $\hat{S}=S-A$ we compute three measures of the size of the \lq\lq error'':
   \begin{itemize}
   \item The $L^{2}$ error. That is to say, the square root of the average
   of $\hat{S}^{2}$.  
   \item The maximum and minimum values of $\hat{S}$.
    \item The \lq\lq maximum and minimum normalised values'' of $\hat{S}$. To motivate this, consider a
    situation where, on  average, $ S $ is very small compared
    with  the full curvature tensor. Then an error
    $\hat{S}$ might appear large compared with the average value of $S$,  while still being small
    compared with the full curvature, and hence it terms of the true geometry
    of the manifold. More generally, we  may have a situation where the
    curvature tensor of the manifold is much larger in some regions than
    others---in other words the natural scales on which to study the manifold
    are different in different regions. To take account of this we define
    the maximum and minimum normalised error to be the maximum and minimum values
    of $\frac{\hat{S}}{\vert {\rm Riem} \vert}$, where $\vert {\rm Riem}\vert$
    is the pointwise norm of the full curvature tensor. 
    
    \end{itemize}

   \subsection{The hexagon}
   
   Recall from the introduction that this is the polygon which corresponds
   to the  Kahler-Einstein metric on the blow-up of the projective
   plane at three general points. In this special case there is another, simpler, approach which one can try---the \lq\lq $T_{K}$ algorithm'' described
in \cite{kn:D1}, and in a separate programme we have verified numerically that this does indeed
converge. But here we will just treat the manifold by our  general procedure---designed for  extremal
metrics. See also the discussion in \cite{kn:Keller}.

We find that we get a useful approximation to the Kahler-Einstein metric
with  the low value $k=6$. The co-efficients of the refined approximation
appear in the figure in Section 2.1, in the table below we record the largest
and smallest of these cefficients and the \lq\lq error measures''.

\

\

$$\begin{tabular}{|l||l|}\hline
Max. Coeff.&$257.2$\\ \hline
Min. Coeff.&$2.062$\\ \hline
$L^{2}$ error & $.0077$\\ \hline
Error maximum &$.021$\\ \hline
Error minimum &$ -.14$ \\ \hline
Normalised error maximum& $.013$\\ \hline
Normalised error minimum& $-.058$ \\ \hline
\end{tabular}$$

\

\

\

 Taking larger values of $k$ we can do much better. There seems to be
 no practical gain in going beyond $k=20$ which yields: 
 
 \
 
 \
 
 \
$$\begin{tabular}{|l||l|}\hline
Max. $a_{\nu}$& $5.52\times 10^{4}$ \\ \hline
Min. $a_{\nu}$& $1.58\times 10^{-4}$\\ \hline
$L^{2}$ error &$1.33\times 10^{-6}$\\ \hline
Error maximum &$1.23\times 10^{-5}$\\ \hline
Error minimum &$-5.01\times 10^{-5}$ \\ \hline
Normalised error maximum&$ 5.8 \times 10^{-6}$\\ \hline
Normalised error minimum& $-2.0\times 10^{-5}$ \\ \hline
\end{tabular}$$

 \
 
 \
 
 In the next table we record data about the size of the components of the
 curvature tensor.
 \
 
  $$ \begin{tabular}{|l||r|}\hline
Max. $\vert {\rm Riem}\vert$ &$2.48$\\ \hline
Min. $\vert {\rm  Riem}\vert$  &$0.61$\\ \hline
Max $K$ & $0$\\ \hline
Min $K$ &$-0.817$\\ \hline 
Max. $\vert w\vert$ & $2.09$\\ \hline
Max. $\vert \rho \vert$ & $2.7\times 10^{-5}$ \\ \hline
\end{tabular}$$
 
 \
 
 \
 
 \
 
 \
 We see that the trace-free Ricci tensor $\rho$ is indeed very small and
 effectively the only non-trivial components
 of the curvature tensor are the \lq\lq Gauss'' and \lq\lq Weyl'' components.
  The Gauss curvature, times $-1/4$, is plotted in Figure 3. The Gauss curvature
 vanishes at the centre of the hexagon and is everywhere negative. (We know
 {\it a priori} that it has to be negative on average, since the Euler characteristic
 of $X_{\bR}$ is $-2$.) The norm of the Weyl component is plotted in Figure
 4. We see that the Weyl component also vanishes at the centre of the hexagon.
 Recall that this Weyl component is a \lq\lq $\phi$-holomorphic'' quartic
 differential on the oriented double cover of $X_{\bR}$. The number of zeros, counted with multiplicity,
 on this double cover is $16$. Since there are $8$ copies
 of $P$ making up this double cover, we can see from symmetry arguments that
 the only possibility is a double zero at the centre of $P$. It is interesting
 to compare the curvature tensor of this manifold with those of the standard
 examples in (5),(6).

\begin{figure}[!h]
  \caption{-1/4 times the Gauss component of the Curvature, for the hexagon}
  \centering
    \includegraphics[scale=0.46]{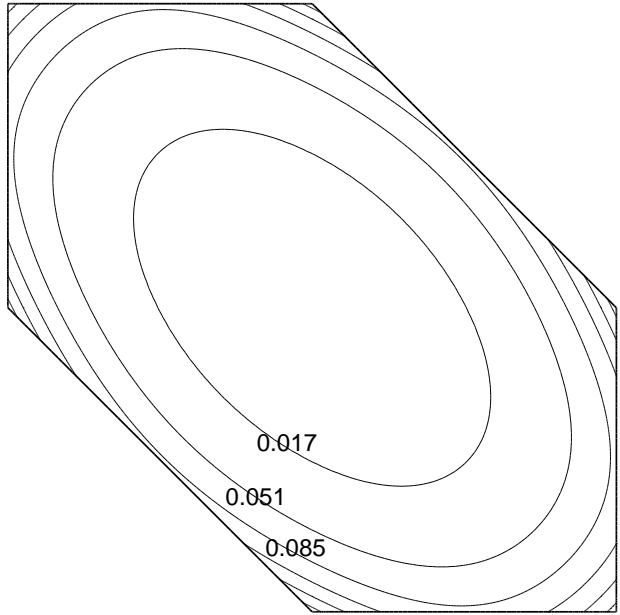}
\end{figure}

\begin{figure}[!h]
  \caption{The norm of the Weyl component of the curvature, for the hexagon}
  \centering
    \includegraphics[scale=0.46]{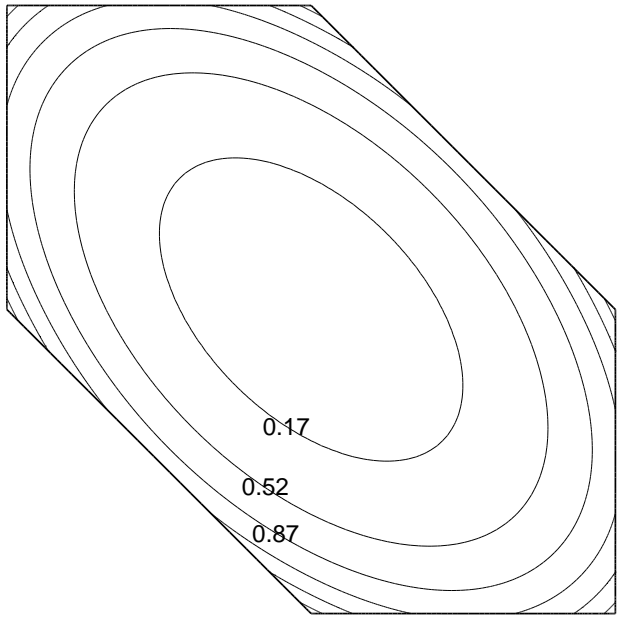}
\end{figure}

     \subsubsection{Geodesics}
     
     Using our numerical solutions for the Kahler-Einstein metric we can
     compute geodesics in this Riemannian manifold. The inner products with
     the two Killing fields  defining the torus action are conserved quantities
     which we refer to as \lq\lq angular momenta''. Fixing the angular momenta,
     the geodesic equations can be regarded as equations for a particle moving
     around the polygon in a potential.  We do not take this any further
     here, and include this discussion only to demonstrate
     that one can effectively study the Riemannian geometry of these solutions.
     One could take the ideas much further. However the solutions  display
     some interesting phenomena. When the angular momenta are reasonably
     large compared with the horizontal component of the velocity the solutions
      behave effectively as an integrable system, essentially confined
     to tori in the phase space. As the angular momenta are reduced the
     dynamics changes its nature and the solutions are more complicated.
     Of course, we might expect the latter from the fact that when the angular
     momenta are zero the equations can be viewed as  the geodesic equations on the
     surface $X_{\bR}$ and since this has negative curvature the system cannot
     be integrable.
     \clearpage
     \begin{figure}[!h]
  \caption{A geodesic in the Kahler-Einstein metric, showing approximately
  integrable behaviour}
  \centering
    \includegraphics[scale=0.7]{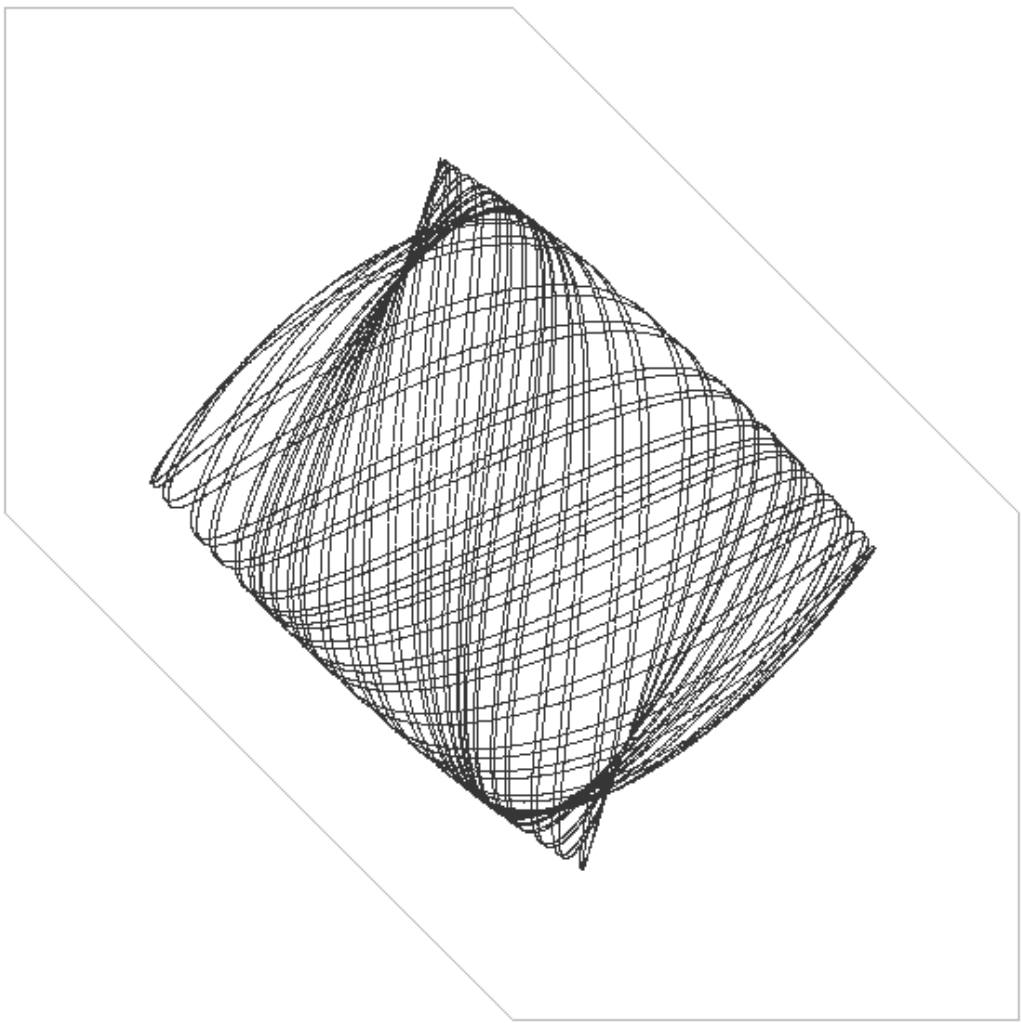}
\end{figure}

     \clearpage
     
     %\newpage
     \subsection{The pentagon}
     
     Here we discuss the pentagon, shown with $k=2$, which corresponds to the blow-up of the
     complex projective plane at two distinct points, with the Kahler class
     equal to $c_{1}$.

    \
    
    \
     
     \

\begin{picture}(150,150)(-70,-20)
\put(0,0){\line(1,0){150}}
\put(0,0){\line(0,1){150}}
\put(150,0){\line(0,1){75}}
\put(0,150){\line(1,0){75}}
\put(75,150){\line(1,-1){75}}
\put(-5,-9){(0,0)}
\put(150,-7){(2,0)}
\put(155,75){(2,1)}
\put(75,5){{\bf 0}}
\put(5,75){{\bf 0}}
\put(140,35){{\bf -1}}
\put(35,141){{\bf -1}}
\put(102,102){{\bf -1}}
\end{picture}
     
     \
     
     \

      The metric we seek is an extremal metric but not of constant scalar
      curvature. The numerical results, regarding speed of convergence and the size of the errors, are very similar to those for the hexagon.

 With $k=20$
we get

 \
 
 \

$$\begin{tabular}{|l||r|}\hline
Max. $a_{\nu}$ &$6.41\times 10^{4}$\\ \hline
Min. $a_{\nu}$ &$3.34\times 10^{-7}$\\ \hline
$L^{2}$ error & $1.02\times 10^{-6}$\\ \hline
Error maximum &.$1.7\times 10^{-5}$\\ \hline
Error minimum & $-6.4\times 10^{-5}$ \\ \hline
Normalised error maximum& $7.4\times 10^{-6}$\\ \hline
Normalised error minimum& $-2.5\times 10^{-5}$\\ \hline
\end{tabular}$$

\

\

The curvature data is given by
  $$  \begin{tabular}{|l||r|}\hline
 Max.$\vert{\rm Riem}\vert$ &$2.57$\\ \hline
 Min.$\vert {\rm Riem}\vert$ &$0.56$\\ \hline
 Max. $K$  & $0.132 $\\ \hline
 Min.$K$ &$-.906$\\ \hline
 Max.$\vert w\vert$ & $2.19$\\ \hline
 Max $\vert \rho\vert $& $.4378$ \\ \hline
\end{tabular}$$

 \
 
 \
     
     The curvature functions are plotted in Figures 6,7,8. (In Figure 8 one
     can see the approximate location of a single zero of the Weyl component,         as
     predicted by topology.)
     
\begin{figure}[!h]
  \caption{The size of the trace-free Ricci tensor, for the pentagon.}
  \centering
    \includegraphics[scale=0.46]{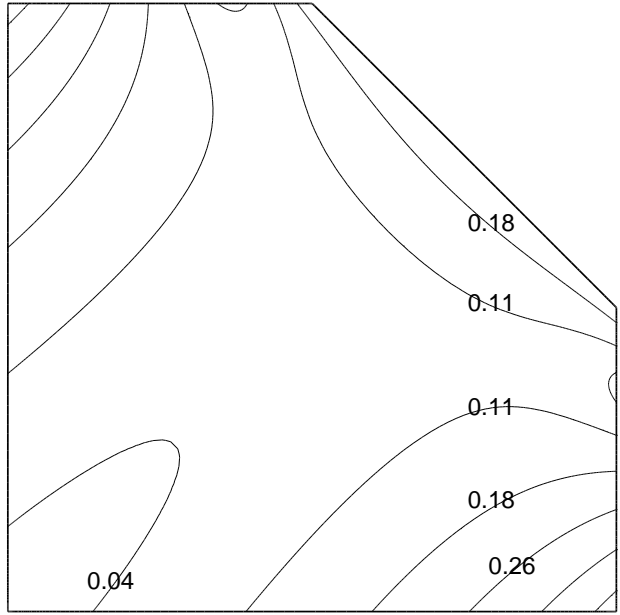}
\end{figure}

\begin{figure}[!h]
  \caption{$-1/4$ times the Gauss curvature, for the pentagon}
  \centering
    \includegraphics[scale=0.46]{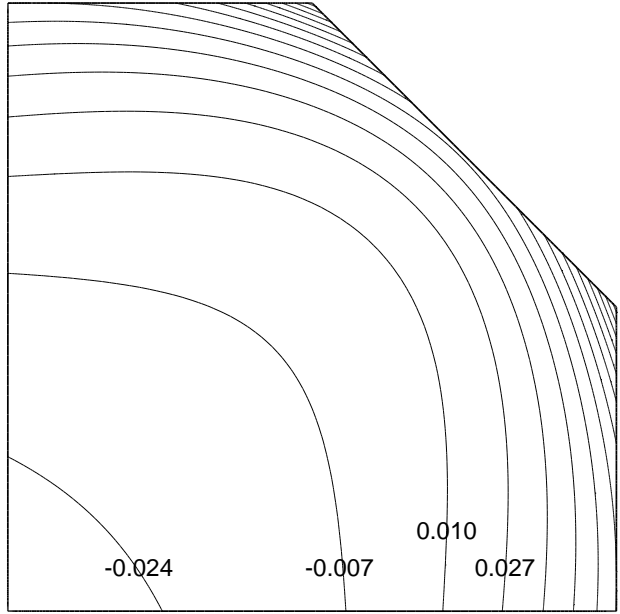}
\end{figure}

\begin{figure}[!h]
  \caption{The size of the \lq\lq Weyl component'', for the pentagon}
  \centering
    \includegraphics[scale=0.46]{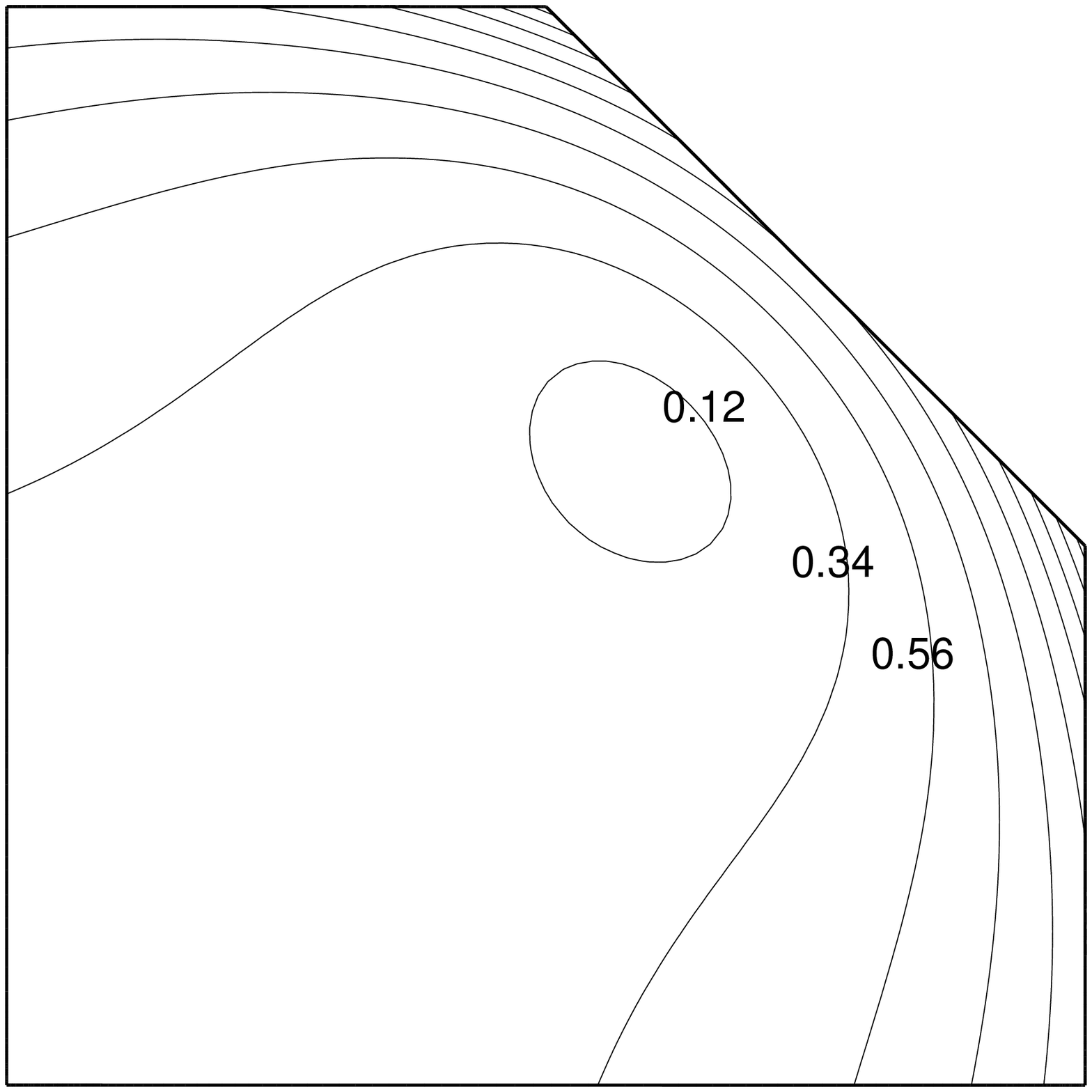}
\end{figure}
     \subsubsection{The conformally-Kahler Einstein metric}
     
     As we mentioned above, the existence of this extremal metric was proved
     recently by Chen, LeBrun and Weber \cite{kn:CBW}. In fact they established
     the existence of a  family of extremal Kahler metrics on this manifold, depending on a real parameter
     $a$ which determines the Kahler class. The class we are
     considering has $a=1$. Chen, LeBrun and Weber show that for a particular value $a_{0}\sim .986\dots$ the Bach tensor vanishes and the extremal metric is conformal to an Einstein
     metric, with conformal factor  the scalar curvature. This extremal metric
     is not directly accessible with our method, which is confined to rational
     values of the parameter with fairly small denominator. However, since
     $a_{0}$ is close to $1$ we can expect that our extremal metric has small
     Bach tensor, leading to a reasonable approximation to the Einstein metric.    We verify this numerically. The norm of the Bach tensor is shown in Figure
   9. It is less than .005 over most of the manifold, and has maximum value
   about $.02$.
   This is comparable with $(1-a_{0})$, as one would expect. It would
   be possible to compute geometric quantities associated to the Einstein metric to
   this degree of accuracy.

     \begin{figure}[!h]
  \caption{The size of the Bach tensor, for the pentagon}
  \centering
    \includegraphics[scale=0.46]{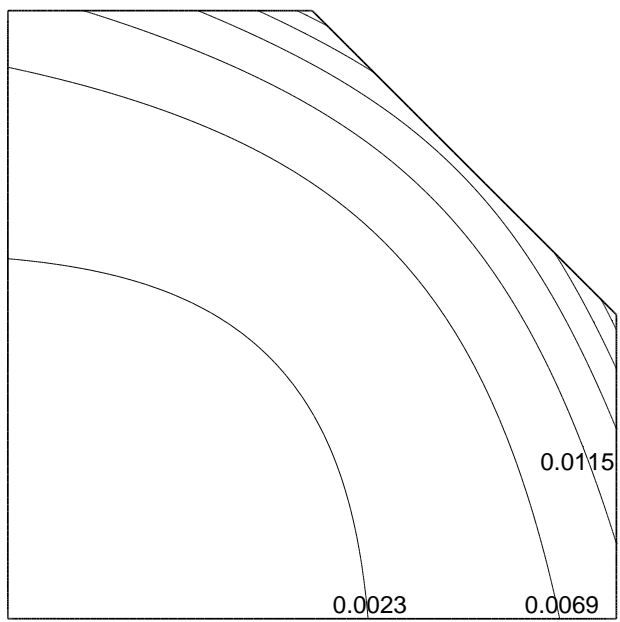}
\end{figure}

     \clearpage
     
     \subsection{The octagon}

     We now move on to the octagon depicted, which corresponds to the projective
     plane blow up at five points, or $S^{2}\times S^{3}$ blown up at four points. The symmetry of the polygon shows immediately that the Futaki invriant
vanishes, so the metric we seek has constant scalar curvature. The existence
of this metric follows from the general result in \cite{kn:D3}. (An argument
of Zhou and Zhu \cite{kn:ZZ} shows that this polygon satisfies the \lq\lq positivity'' criterion in
the general existence theorem.)

\

\

  \
 
  \
   
\begin{picture}(180,180)(-70,-30)
\put(0,60){\line(0,1){60}}
\put(0,120){\line(1,1){60}}
\put(60,0){\line(1,0){60}}
\put(120,0){\line(1,1){60}}
\put(180,60){\line(0,1){60}}
\put(60,180){\line(1,0){60}}
\put(180,120){\line(-1,1){60}}
\put(60,0){\line(-1,1){60}}
\put(32,32){{\bf -1}}
\put(88,3){{\bf -2}}
\put(5,88){{\bf -2}}
\put(142,35){{\bf -1}}
\put(142,142){{\bf -1}}
\put(168,90){{\bf -2}}
\put(88,172){{\bf -2}}
\put(32,142){{\bf -1}}
\put(60,-8){(2,0)}
\put(120,-7){(4,0)}
\put(180,55){(6,2)}
\end{picture}

\

\

\

We find that we need to take $k$ significantly larger than in the previous
two cases, to get a reasonable approximate solution. Taking $k=18$ we obtain
     
$$\begin{tabular}{|l||r|}\hline
Max. $a_{\nu}$ &$77349.5$\\ \hline
Min. $a_{\nu}$ &$.00051733$\\ \hline
$L^{2}$ error & $.00522$\\ \hline
Error maximum &$.0334$\\ \hline
Error minimum & $-.0918$\\ \hline
Normalised error maximum&$.011$\\ \hline
Normalised error minimum& $-.0135$\\ \hline
\end{tabular}$$
     
     \
     
     \

     When $k=36$  we get
     
\

\

 $$    \begin{tabular}{|l||r|}\hline
Max. $a_{\nu}$.&$2.16\times 10^{8}$\\ \hline
Min. $a_{\nu}$.&$9.94\times 10^{-10}$\\ \hline
$L^{2}$ error & $1.07\times 10^{-4}$\\ \hline
Error maximum &$1.32\times 10^{-3}$\\ \hline
Error minimum & $-3.48\times10^{-3}$\\ \hline
Normalised error maximum &$2.69\times 10^{-4}$\\ \hline
Normalised error minimum& $-5.07\times 10^{-4}$ \\ \hline
\end{tabular}$$
    
    \
    
    \

     The reason why we need to make $k$ larger than for the previous two
     polygons is clear when we look at the
     size of the curvature tensor:
     
     \
     
     \
     
   $$     \begin{tabular}{|l||r|}\hline
 Max. $\vert {\rm Riem}\vert$ &$6.86$\\ \hline
 Min.$\vert {\rm Riem}\vert$ &$.644$\\ \hline
 Max. $K$ & $0$\\ \hline
Min.$K$ &$-2.60$\\ \hline
 Max.$\vert w\vert$ & $5.87$\\ \hline
Max $\vert \rho \vert$ & $1.35$ \\ \hline
\end{tabular}$$

\

\

The maximum curvature in this case is nearly 3 times that for the hexagon
and pentagon. Now $\sqrt{k}$ behaves as a length-scale parameter in the asymptototic
theory, and curvature scales as ${\rm length}^{-2}$. So we would expect that
for the octagon we should take $k$ about 3 times as large as for the pentagon
or hexagon, to obtain comparable accuracy. The figures above are consistent
with this. The size of the full curvature tensor and of the trace-free Ricci
component are shown in Figures 10,11. Notice that the curvature is much larger near the boundary of the polygon. Moreover it is larger  near
     edges corresponding to  the $-2$ curves than  the $-1$ curves. The Chern-Weil
     formula (7) shows that as the number of edges of the polygon increases
     the full curvature tensor grow, on average, compared to the scalar curvature.
     Also, the self-intersection number of a boundary sphere can be
     written as an integral of an expression involving the curvature over
     the corresponding edge, so the curvature must become relatively large
     there. But the actual concentration of the curvature is more dramatic
     than one might have expected on such general grounds.
     
     Figure 12 shows the size of the \lq\lq error'' $S-A$ over the polygon.
     As we would expect the error is much larger near the boundary, where
     the curvature is large (which is reflected by the normalised errors
     in the tables above). Notice the intricate oscillatory pattern in the
     function, which one can see is related to the lattice spacing.

          \begin{figure}[!h]
  \caption{$\vert {\rm Riem}\vert$ for the octagon}
  \centering
    \includegraphics[scale=0.8]{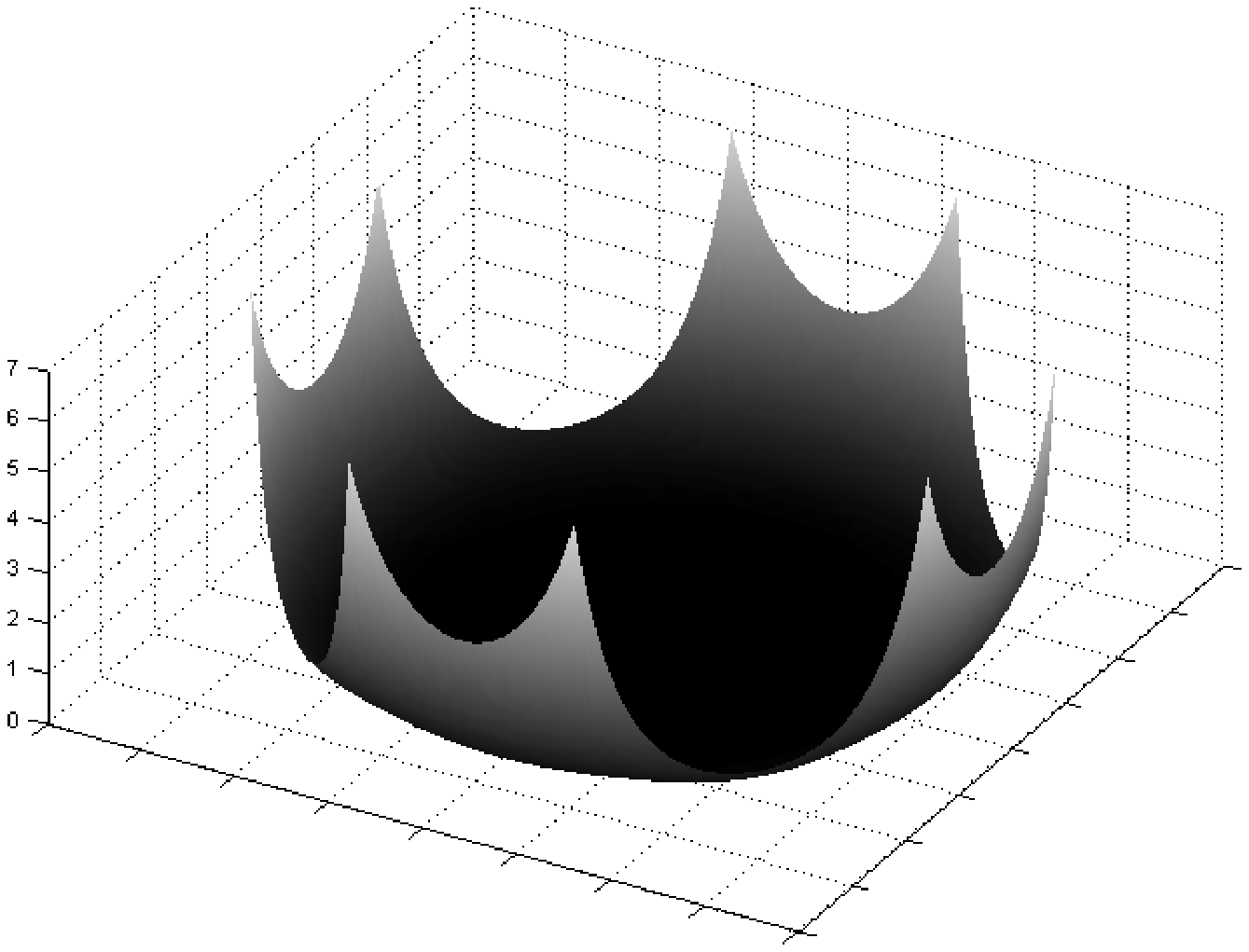}
\end{figure}

\begin{figure}[!h]
  \caption{The size of the trace-free Ricci component}
  \centering
    \includegraphics[scale=0.8]{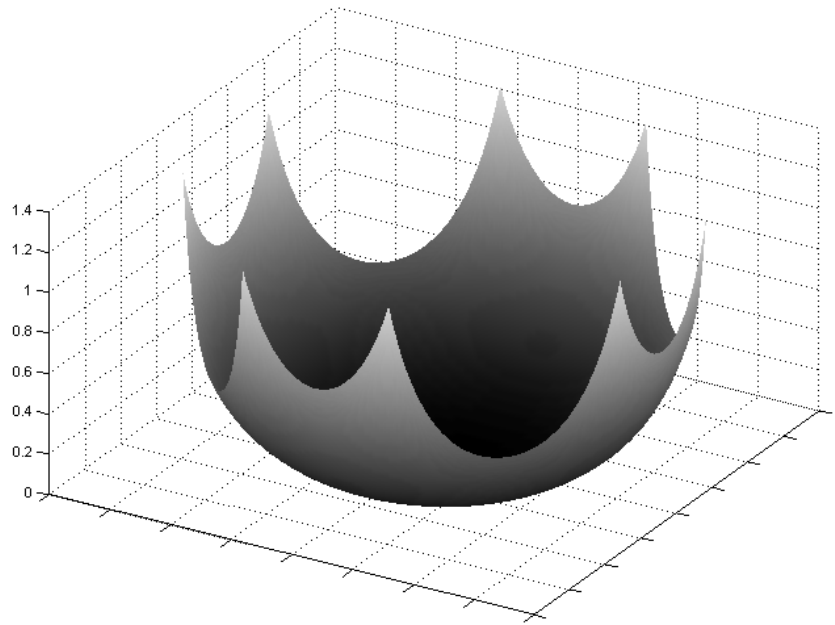}
\end{figure}

\begin{figure}[!h]
  \caption{The size of the error $\hat{S}$}
  \centering
    \includegraphics[scale=0.8]{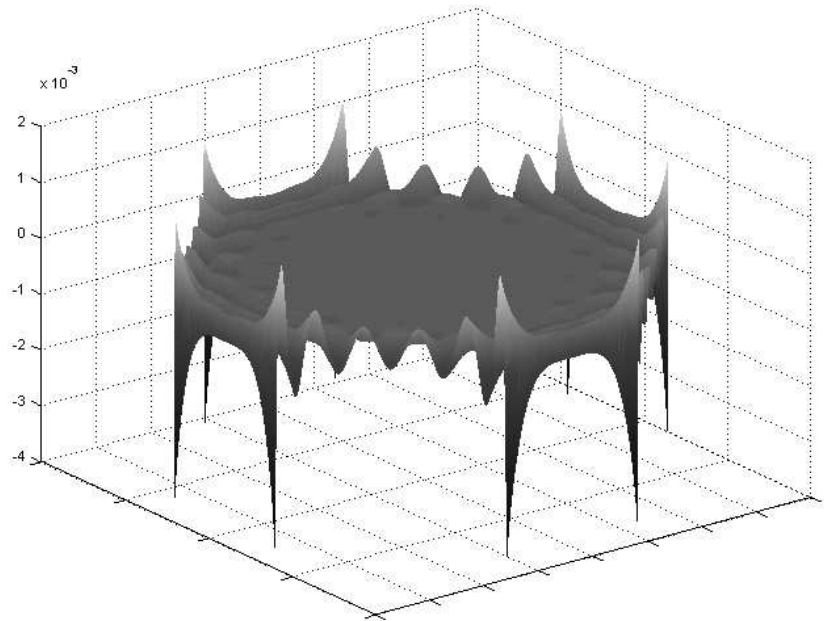}
\end{figure}
\clearpage
     
     \subsection{The heptagon}

     This is the polygon as shown (with $k=5$). It corresponds to a manifold  we get by blowing up two of the fixed points in the pentagon manifold. The
 self-intersection numbers are as shown, and there is one $-3$ curve.

     \
     
     \
      
      \

      \begin{picture}(200,200)(-70,-30)
\put(0,0){\line(1,0){200}}
\put(0,0){\line(0,1){200}}
\put(200,0){\line(0,1){40}}
\put(200,40){\line(-1,2){40}}
\put(160,120){\line(-1,1){40}}
\put(120,160){\line(-2,1){80}}
\put(40,200){\line(-1,0){40}}
\put(5,90){{\bf 0}}
\put(90,5){{\bf 0}}
\put(130,130){{\bf -3}}
\put(190,20){{\bf -2}}
\put(20,190){{\bf -2}}
\put(175,65){{\bf -1}}
\put(65,175){{\bf -1}}
\put(0,-9){(0,0)}
\put(200,-9){(5,0)}
\put(205,40){(5,1)}
\put(165,120){(4,3)}
\put(120,165){(3,4)}
\end{picture}

\

\

The metric we seek is an extremal metric, not of constant scalar curvature. The results follow
      a similar pattern to the previous case, only more so. We take $k=45$
      (which is about as large as seems reasonable with the programme), and
      we obtain:

      $$   \begin{tabular}{|l||r|}\hline
Max. $a_{\nu}$&$3.88\times 10^{9}$\\ \hline
Min. $a_{\nu}$&$2.44\times 10^{-21}$\\ \hline
$L^{2}$ error & $2.6\times 10^{-3}$\\ \hline
Error maximum &$ .085$\\ \hline
Error minimum &$-.23$ \\ \hline
Normalised error maximum& .015\\ \hline
Normalised error minimum& -.016 \\ \hline
\end{tabular}$$

      The error is concentrated near the boundary and particularly near the
      $-3$ curve, which is where the curvature
      is largest. The modulus $.23$ of the minimum value looks large, and one
      might think
      that the approximation is not very satisfactory, but the normalised
      value is much smaller. Thus we believe that this does model an exact
      solution quite well. The maximum value of the curvature is  nearly
      6 times that for the hexagon, so we should expect to achieve an accuracy
      similar to that for the hexagon with $k\sim 45/6$, and the results are
      in line with this. Notice that the curvature function in Figure 13
      is almost constant, at a value about $0.7$,  over most of the polygon.

  $$    \begin{tabular}{|l||r|}\hline
 Max.$\vert{\rm Riem} \vert$ &$13.9$\\ \hline 
 Min.$\vert {\rm Riem}\vert$ &$.537$\\ \hline
 Max. $K$ & $.249 $\\ \hline
 Min.$K$  &$-5.48$\\ \hline
 Max. $\vert w\vert$  & $11.7$\\ \hline
Max. $\vert \rho \vert$ & $ 4.27$ \\ \hline
\end{tabular}$$

\

\

\begin{figure}[!h]
  \caption{$\vert {\rm Riem}\vert$ for the heptagon}
  \centering
    \includegraphics[scale=0.8]{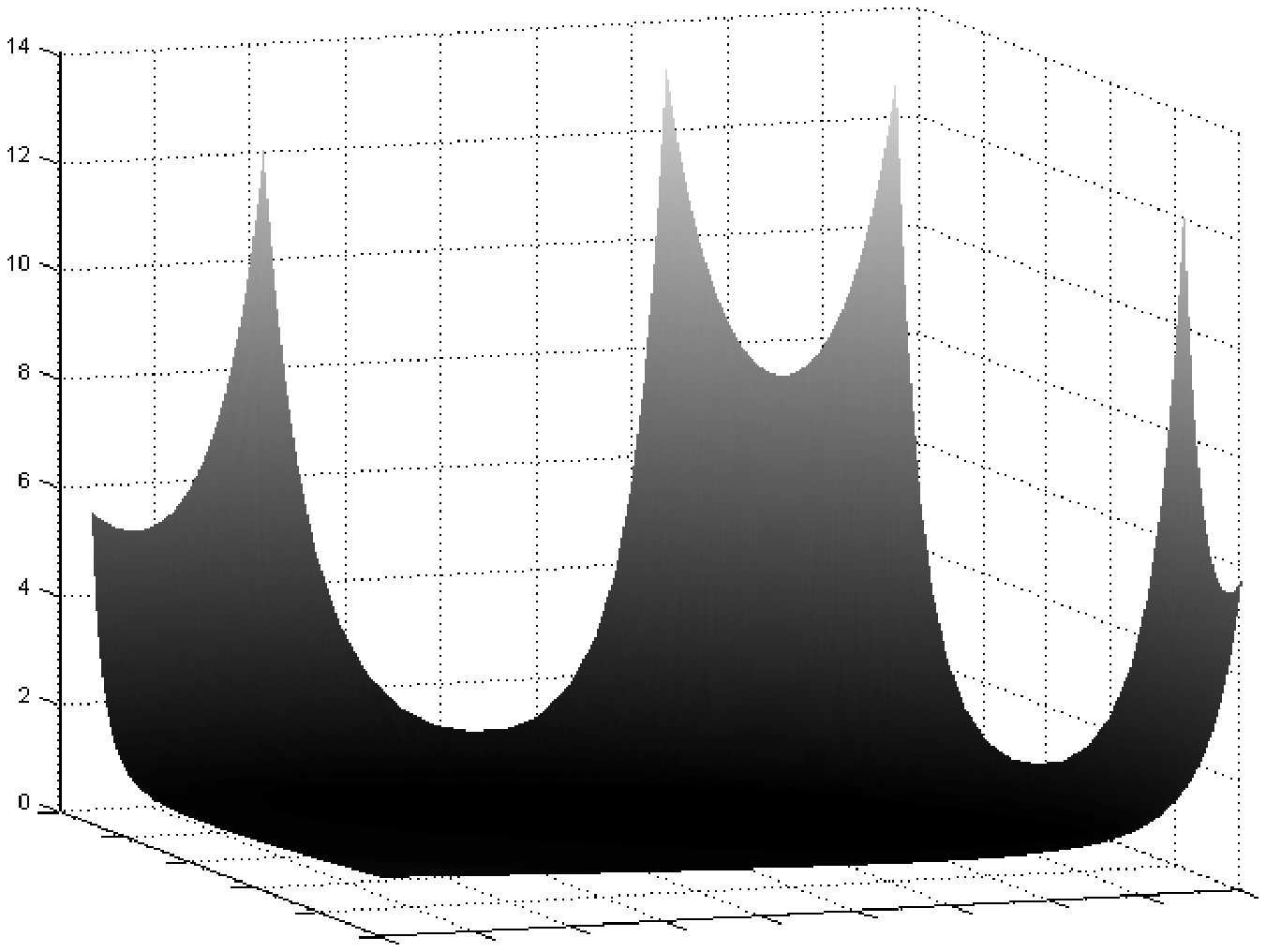}
\end{figure}
   
%\begin{figure}[!h]
%  \caption{heptagon bach}
%  \centering
%    \includegraphics[scale=0.8]{heptagon_bach}
%\end{figure}

\begin{figure}[!h]
  \caption{The error $\hat{S}$, for the heptagon}
  \centering
    \includegraphics[scale=0.8]{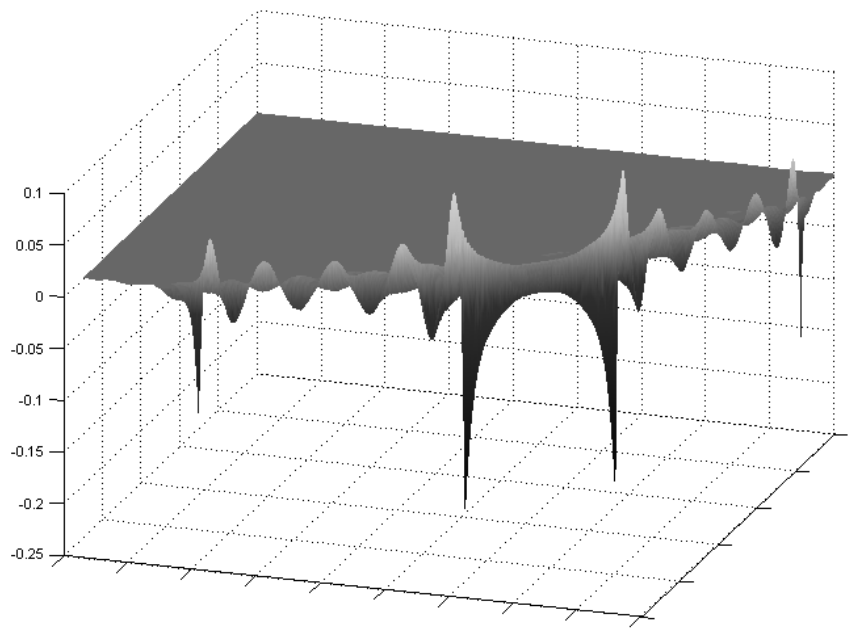}
\end{figure}

     \clearpage
     \section{Conclusions}

     \begin{itemize}
     \item {\it Rate of convergence}
     
     As we have stated in Section 3, one can expect that the refined approximations
     converge to the extremal metric faster than any inverse power. But,
     even if proved rigorously, this is only an asymptotic statement and
     does not give any information about particular values of $k$. To examine
     this experimentally, we consider the $L^{2}$ error as a measure of
     the distance to the true solution, working with   the hexagon,
     over a range of different values of $k$:
     
     $$\begin{tabular}{l|r|r|r|r|r|r|r} \\ \hline
    $ k$&4&6&8&10&12&16&20\\ \hline
     Error &.019&.0077&.0012&.00052&.00010&$8.0\times 10^{-6}$&$ 1.3\times 10^{-6}$\\ \hline
     \end{tabular}$$
     
     The data fits reasonably well with an {\it exponential} decrease $10^{-(k/4)-1}$.

     \item{Speed and size of the calculation}
     For a fixed \lq\lq shape'' of polygon the number of elementary calculations required in each step of the iteration
     is proportional to $N^{2} k^{2}$, since we have to evaluate $O(k^{2})$
     functions, and their derivatives at $O(N^{2})$ points. In the more general
     case of a $d$-dimensional toric variety this would be $N^{d}k^{d}$.
     As $k$ increases
     we need to increase $N$, since the functions we are integrating become
     increasingly localised. Simple asymptotics suggests that one should
     take
     $N=O(\sqrt{k})$, so the number of elementary calculations is $O(k^{3d/2})$.
     However, at a given point, most of the functions are extremely small, so
     one can hope to reduce the calculation by ignoring them. This idea was
     implemented to some extent in our programme, which \lq\lq truncates''
     the sums to ignore the very small terms. This does not affect the accuracy
     noticeably but the gain we have achieved so far is only modest. With
     a more sophisticated approach, one could hope to make substantial gains,
     particularly when $k$ is large, or in higher dimensional problems.
     Simple asymptotics suggests that one might be able to reduce the calculation
     at each point to $O(k^{d/2})$, giving $O(k^{d})$ for each iteration step.

     The convergence is also slower as $k$ increases. One expects that the
     number of steps required to effectively reach a fixed point will be
     $O(k^{2})$. (For example, in Figure 2 the number of steps to reach the
     marked \lq\lq corner'' in the graph for the pentagon, with $k=20$, is
     about $500$, while for the heptagon, with $k=45$, it takes about $2200$
     steps: then  $2200/500 \sim (45/20)^{2}$.) So in sum we can hope that the length of time needed to run the whole
     procedure, in dimension $d$, will be $O(k^{2+3d})$ without the truncation of the sums described above,
     and perhaps $O(k^{2+2d})$, with this truncation.

     \item{Scope for improvement}
     
     The most striking feature of our results is the extent to which the
     curvature of the solutions becomes concentrated around the boundary
     of the polygon or, in terms of the $4$-manifolds, around the exceptional
     curves. This is particularly true of the curves with large negative
     self-intersection. The obvious way to handle this efficiently would
     be to adopt a \lq\lq multiscale'' technique. This idea can be related
     to
     the \lq\lq Veronese embedding'' in algebraic geometry, as discussed
     in \cite{kn:D2}. Suppose we have an algebraic metric $\omega$ defined by an array of coefficients
     $a_{\nu}$, for a certain lattice-size $k_{0}$. The Kahler potential is $\phi= \log (\sum a_{\nu} e^{\nu t})$. Then we can write
$$  \phi= \frac{1}{2} \log\left( (\sum a_{\nu} e^{\nu t})^{2}\right).$$
We can expand out the sum to write
$$  \left(\sum a_{\nu} e^{\nu t}\right)^{2}=   \sum B_{\mu} e^{\mu t}, $$
where 
$$  B_{\mu}= \sum_{\nu_{1}+\nu_{2}= \mu}  a_{\nu_{1}} a_{\nu_{2}} .$$
 This means that, after rescaling, we can
obtain the {\it same} metric $\omega$ using the array $B$ but with lattice
scale $2k_{0}$. Now suppose we have chosen a subset $G\subset 2k_{0} P \cap
\bZ^{2}$. Then we can consider varying our array to
$$   \tilde{B}_{\mu}= B_{\mu} +c_{\mu}, $$
where $c_{\mu}$ is zero if $\mu$ is not in $G$. The potential function
is
$$   \log\left( (\sum a_{\nu} e^{\nu t})^{2} + \sum c_{\mu} e^{\mu t}\right).
$$
If the number of elements of $G$ is relatively small then we can calculate
this function and its derivatives relatively quickly. But if we choose $G$
appropriately we can hope to get a much more accurate approximation to the
true metric over the high curvature regions. For example in the case of the
heptagon we might take $k_{0}=30$, which should suffice to describe the metric
over the interior of the polygon, and a set $G$ consisting of lattice points
near to the boundary,  so that we would effectively
be approximating the metrics in that region as though we were working with
$k=60$. 

Of course we can imagine iterating this construction to a yet larger lattice
scale $4k_{0}$, and so on. Also, there is nothing special about the factor
$2$. But what is needed to develop this approach is to find some reliable
algorithms for updating the different coefficients.

\end{itemize}

%+Bibliography

%-Bibliography

\end{document}